\documentclass[12pt,sort&compress]{elsarticle}

\usepackage{setspace}
\usepackage[top=1in,bottom=1in,left=1in,right=1in]{geometry}

\usepackage{amssymb}
\usepackage{mathrsfs}
\usepackage{amsmath}
\usepackage{amsfonts}
\usepackage{amsthm}
\usepackage{color}
\usepackage{subcaption}
\usepackage[skip=4pt]{caption}
\usepackage{float}
\usepackage{url}
\usepackage{multicol}

\usepackage[pdfborder={0 0 0},colorlinks,allcolors=blue]{hyperref}
\usepackage{newtxtext, newtxmath}
\usepackage{arydshln}
\usepackage{enumitem}
\usepackage{gensymb}

\usepackage{nth}
\usepackage{xfrac}

\usepackage{comment}

\usepackage{algorithm}
\usepackage{algpseudocode}
\usepackage[export]{adjustbox}
\usepackage{multirow}
\usepackage{booktabs}

\theoremstyle{definition}

\graphicspath{{figures_arxiv/}}

\begin{document}
\begin{frontmatter}

\title{Accurate wall shear stress in immersed flow analysis with application to point cloud-based CFD}

\author[inst1]{Monu~Jaiswal}
\author[inst1]{Ming-Chen~Hsu\corref{cor}}
\ead{jmchsu@iastate.edu}

\cortext[cor]{Corresponding author}

\address[inst1]{Department of Mechanical Engineering, Iowa State University, 2043 Black Engineering, Ames, Iowa 50011, USA}

\begin{abstract}
Point cloud-based CFD enables flow analysis directly on discrete points obtained from 3D scanning and medical imaging, bypassing surface reconstruction, geometry cleanup, and boundary-fitted mesh generation. Derived from immersogeometric analysis, the method immerses the point cloud in a background mesh and enforces no-slip conditions on discrete points through a Nitsche-based weak boundary condition (BC). The framework delivers accurate velocity fields, pressure distribution, and integrated loads; however, accurate prediction of the local wall shear stress (WSS) has remained a critical challenge. The geometry intersects the background mesh arbitrarily, producing cut elements that lack the regularity required for consistent gradient evaluation. The issue is compounded by the stabilization term of the weak BC, whose parameter estimation in the symmetric Nitsche formulation is dependent on the cut configuration and affects the variationally consistent definition of traction from which the WSS is computed. In this work, we propose a new method to obtain accurate wall shear stress in immersed flow analysis with application to point cloud-based CFD, using a non-symmetric Nitsche's formulation with near-wall modeling and a patch-based stress recovery approach with traction compatibility. The method is validated on canonical benchmarks and applied to turbulent flow past a sphere and to a patient-specific aorta, showcasing excellent agreement with reference results.

\end{abstract}


\end{frontmatter}

\tableofcontents

\section{Introduction}
\label{sec:intro}
Computational Fluid Dynamics (CFD) has become an indispensable tool across engineering disciplines~\cite{Schwarz23Beyon, Hughes20Compu, Bazilevs19Compu2, Korobenko19Compu1, Silvestri21CFD, Updegrove17SimVa} for analyzing complex flow behavior and extracting quantities of interest, such as lift, drag, and wall shear stress (WSS), for crucial design decisions and diagnostics. However, traditional CFD workflows require a boundary-fitted mesh, which entails surface reconstruction, geometry cleanup, and de-featuring to obtain a watertight model, followed by surface and volume mesh generation under strict mesh quality criteria~\cite{Dawes01Reduc, Spalart16Onrol}. These preprocessing steps are tedious, labor-intensive, and frequently dominate the turnaround time of the entire engineering analysis process~\cite{Wolfe09Immer, Bazilevs10Isoge}. Immersed boundary methods (IBM) circumvent boundary-fitted meshing altogether by immersing the geometry in a non-conforming background mesh~\cite{Peskin02immer, Mittal05Immer, Griffith09Simul}. Within finite element-based immersed methods~\cite{Duster08finit, Burman15CutFE, Massing14stabi, Schott14new, Main18shift2}, immersogeometric analysis (IMGA)~\cite{Kamensky15immer} offers a variationally consistent approach in which the flow is discretized using a variational multiscale (VMS) formulation~\cite{Bazilevs07Varia, Hughes00Large} and no-slip Dirichlet conditions are enforced on the immersed surface with Nitsche-based weak boundary conditions (BC)~\cite{Bazilevs07Weak1, Jaiswal26Weak, Hsu26revie}. IMGA accommodates diverse geometric representations, including non-uniform rational B-splines (NURBS)~\cite{Kamensky15immer}, boundary representation (B-Rep)~\cite{Hsu16Direc}, polygonal meshes~\cite{Xu16tetra}, and has recently been extended to point clouds~\cite{Balu23Direc}. 

The point cloud-based CFD framework~\cite{Balu23Direc} enables direct analysis by immersing the point cloud in a tetrahedral background mesh and imposing the weak BC directly on the discrete points, without any surface reconstruction. The framework not only eliminates geometry cleanup for dirty CAD models, which can be sampled directly into point clouds, but also establishes a seamless scan-to-CFD pipeline for real-world geometries where no CAD model exists. For example, it has been applied to scanned point clouds of in-use civil structures acquired with photogrammetry~\cite{Wang23Photo} and to patient-specific point clouds extracted from medical images through deep learning-based auto-segmentation models~\cite{Corpuz25Direc, Corpuz26Multi}. To improve robustness with scanned point clouds, which exhibit noise, gaps, and non-uniform sampling, the framework was further combined with ghost penalty stabilization for small cut elements and a mesh-driven resampling strategy~\cite{Jaiswal24Mesh}. Across these studies, point cloud-based CFD has been extensively validated and applied to complex geometries that were previously impractical to analyze. The framework delivers field solutions such as velocity and pressure, and integrated loads such as drag and lift, with accuracy comparable to boundary-fitted simulations at similar degrees of freedom. Despite these advantages, obtaining accurate and smooth wall shear stress remains a challenge for point cloud-based CFD, and for immersed methods in general~\cite{Zhang22impro, Capatina21Flux}.

Wall shear stress quantifies the tangential traction exerted by the fluid on the solid boundary. In aerodynamic applications, accurately resolving WSS is essential for predicting skin friction, detecting flow separation, and modeling convective heat transfer~\cite{Pope00Turbu, Bergman18Funda}. WSS is equally vital in cardiovascular flow, where abnormal shear stress acts as an indicator for endothelial dysfunction, atherosclerosis, and aneurysm progression~\cite{Kumar18Lowco, Zhou23Walls, Shojima04Magni, Boyd16Lowwa}. Despite its importance, robust WSS estimation remains challenging because it is a gradient field computed from the primal velocity solution, making it highly sensitive to mesh quality~\cite{Brunatova25Onnum, Berzins99Meshq}. For the boundary-fitted method, accurate WSS estimation relies on strict boundary-layer refinement, with wall-aligned elements of controlled wall-normal size, ensuring smooth, resolved velocity gradients at the wall~\cite{Takizawa10Walls}. In the immersed setting, no such structure exists; the geometry intersects the background mesh arbitrarily, producing cut elements of uncontrolled effective size and shape.  This unstructured, non-aligned discretization lacks the regularity required for consistent gradient evaluation, and the result is oscillatory and sub-optimal WSS, even when the primal fields and integrated forces are accurate.  

Moreover, in the IMGA approach, the WSS is extracted from the variationally consistent definition of boundary traction~\cite{Bazilevs10Large, vanBrummelen12Flux}: $ - \boldsymbol{\sigma}^h \mathbf{n}- \rho \left\{ \mathbf{u}^h \cdot \mathbf{n}\right\}_{-} \left(\mathbf{u}^h - \mathbf{g} \right) + \tau^B \left( \mathbf{u}^h - \mathbf{g} \right)$, where $\boldsymbol{\sigma}^h$ is the Cauchy stress tensor, $\mathbf{n}$ is the outward normal, $\rho$ is the density, $\mathbf{u}^h$ is the velocity, $\mathbf{g}$ is the prescribed Dirichlet data, $\left\{\cdot\right\}_{-}$ denotes the negative part of its argument, and $\tau^B$ is the stabilization parameter of Nitsche's method. The first term is the Cauchy traction whose viscous part involves the velocity gradient and therefore inherits the evaluation difficulties described above. The second term acts only on the inflow portion of the boundary and scales with the BC residual, making it the least dominant. The third term, however, is the stabilization/penalty term of Nitsche's method, which is independently compromised. In the symmetric Nitsche's method, coercivity imposes a lower bound on the parameter $\tau^B$ derived from the discrete trace inverse inequality~\cite{Evans13Expli}. For a boundary-fitted mesh, this bound takes a simple form $\tau^B = C\mu/h$, where $h$ is the wall-normal element size~\cite{Bazilevs07Weak1}. Therefore, using a boundary layer mesh with uniform wall-normal spacing yields a stabilization contribution that is uniform along the wall. However, determining the lower bound of $\tau^B$ for IMGA is not straightforward and involves solving the element-local eigenvalue problem~\cite{Embar10Impos} which is highly sensitive to cut configurations~\cite{Schillinger16nonsy}. Sliver cuts and small volume fractions inflate the stabilization parameter by orders of magnitude element-by-element. As a result, even a small BC residual, once multiplied by a parameter that varies significantly between elements, contributes a spatially discontinuous and noisy field to the traction. The alternative of a large constant parameter satisfying the bounds everywhere (for example, $\tau^B = 10^3$ \cite{Xu16tetra}) does not offer an advantage either; it uniformly amplifies the non-physical BC residual, so that the traction is dominated by penalty noise rather than physical stress. Both significant terms of the consistent traction are compromised: the gradient by the geometry of the cut and the penalty by the parameter required for stability. 

Extracting accurate WSS and skin friction from non-conforming and immersed methods remains an active and challenging subject of study. In the finite volume community, standard sharp-interface IBMs impose boundary conditions by reconstructing near-wall velocity via interpolation from surrounding fluid nodes, and skin friction is then extracted from near-wall stress balance or explicit finite-difference approximations of the wall-normal gradient~\cite{Tseng03ghost, Mittal08versa, Roman09simpl}. Within finite element frameworks, evaluating an accurate boundary stress on cut elements requires higher-degree bases~\cite{Duster08finit} or local boundary-fitted meshes near the immersed surface~\cite{Massing15Nitsc}. Stress recovery techniques~\cite{Guo25Recov}, particularly the superconvergent patch recovery (SPR) of Zienkiewicz and Zhu~\cite{Zienkiewicz92super} with equilibrium and boundary enhancements~\cite{Blacker94Super, Wiberg94Super, Wiberg95Impro}, are often used to reconstruct an improved gradient field from interior solutions. These techniques have been extended to immersed methods, such as locally conservative flux recovery for cut finite element solutions of elliptic problems~\cite{Capatina21Flux}, stress improvement procedures for enriched and unfitted formulations in solid mechanics~\cite{Zhang22impro}, and SPR with equilibrium and traction-continuity constraints in immersed Cartesian-grid contact problems~\cite{Navarro20Super}. Notably, Tur et al.~\cite{Tur15Stabi} employed the stress recovery technique not merely as a post-processing step, but within the boundary condition formulation itself to stabilize the Dirichlet imposition on the immersed meshes. However, these recovery efforts concern elliptic and solid-mechanics problems and only target the recovery of the gradient term. To the best of our knowledge, no existing method addresses the combined effect of the gradient and penalty terms of Nitsche's method on computing variationally consistent WSS in immersed methods. 

In this work, we propose a method for computing accurate and smooth WSS for point cloud-based CFD by coupling the non-symmetric Nitsche formulation~\cite{Schillinger16nonsy, Burman12penal, Boiveau16penal} for weak BC with a patch-based WSS recovery technique inspired by SPR~\cite{Zienkiewicz92super, Blacker94Super}. In contrast to the symmetric formulation, the non-symmetric Nitsche's method is coercive for any penalty parameter $\tau^B \geq 0$, thereby eliminating the strict lower bound required for numerical stability~\cite{Schillinger16nonsy}. However, while stability is guaranteed, introducing a strictly positive penalty ($\tau^B > 0$) remains essential to ensure accuracy for complex problems~\cite{Guo17param, Heimann13unfit}. To appropriately determine this parameter, we employ a near-wall model based on Spalding's formula for the law of the wall~\cite{Spalding61singl}. Using near-wall modeling for weak BC parameter estimation has shown improved accuracy for symmetric Nitsche's method~\cite{Golshan15Large}. Next, we introduce a patch-based WSS recovery strategy that combines the interior velocity solution with the variationally consistent traction obtained from the weak BC. To prevent non-physical diffusion of the steep near-wall velocity gradient, as is evident in high-Reynolds number flows, we constrain the recovery using a restrictive, element-based patch construction. Together with our point cloud-based CFD, these advancements form a robust framework to perform direct flow analysis on point cloud geometries while yielding accurate field solutions, integrated forces, and local WSS distribution.

The remainder of this paper is organized as follows. Section~\ref{sec:framework} presents the computational framework, including the IMGA formulation for incompressible flow with ghost penalty stabilization, non-symmetric Nitsche's method for weak BC, and WSS recovery procedure. Section~\ref{sec:validation} validates the method on two canonical benchmarks: Hagen--Poiseuille pipe flow at $Re=50$ compared with the exact solution and flow around a sphere at $Re=100$ compared against a boundary-fitted result. Section~\ref{sec:applications} presents the framework for turbulent flow past a sphere at $Re=3700$ and a patient-specific aorta, demonstrating excellent agreement with reference results. Section~\ref{sec:conclusion} draws conclusions and outlines directions for future research.

\section{Framework}  
\label{sec:framework}
This section presents the complete framework for direct flow analysis on point cloud geometry. First, we give a brief overview of the immersogeometric formulation for incompressible flow~\cite{Kamensky15immer, Xu16tetra} and ghost penalty stabilization for small cut elements~\cite{Burman10Ghost}. Next, we describe the two components developed for accurate WSS calculation: non-symmetric Nitsche-based weak BC and patch-based WSS recovery. Finally, a cut quadrature rule constructed directly from point clouds is presented for performing accurate integration on the background mesh. 

\subsection{Immersogeometric formulation for incompressible flows}
\label{SubSec:Imcompressible}

\begin{figure}[!t]\centering
  \begin{subfigure}{0.3\textwidth}\centering
    \includegraphics[width=\textwidth]{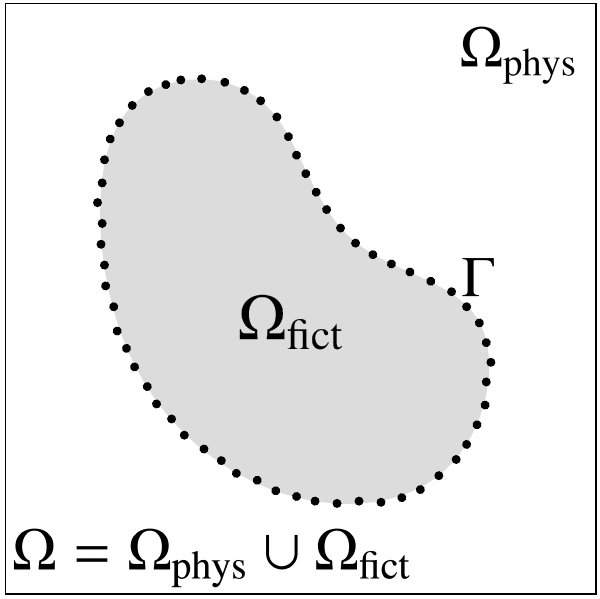}
    \caption{}
    \label{fig:IMGA_domains}
  \end{subfigure}
  \hspace{0.01\textwidth}
  \begin{subfigure}{0.3\textwidth}\centering
    \includegraphics[width=\textwidth]{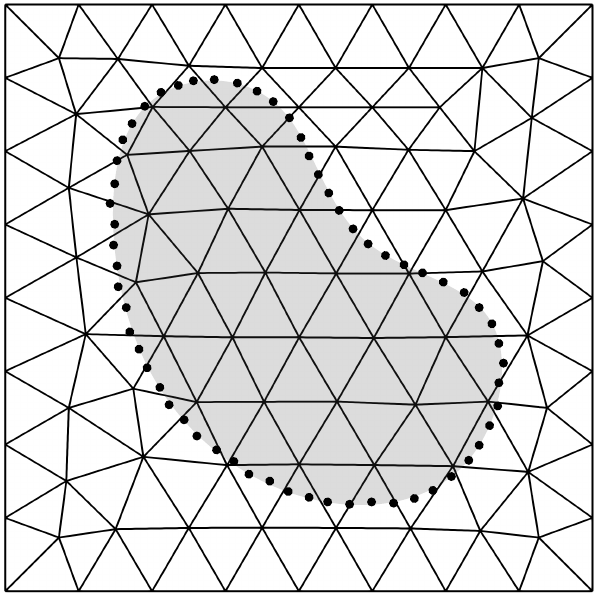}
    \caption{}
    \label{fig:IMGA_mesh}
  \end{subfigure}
  \hspace{0.01\textwidth}
  \begin{subfigure}{0.3\textwidth}\centering
    \includegraphics[width=\textwidth]{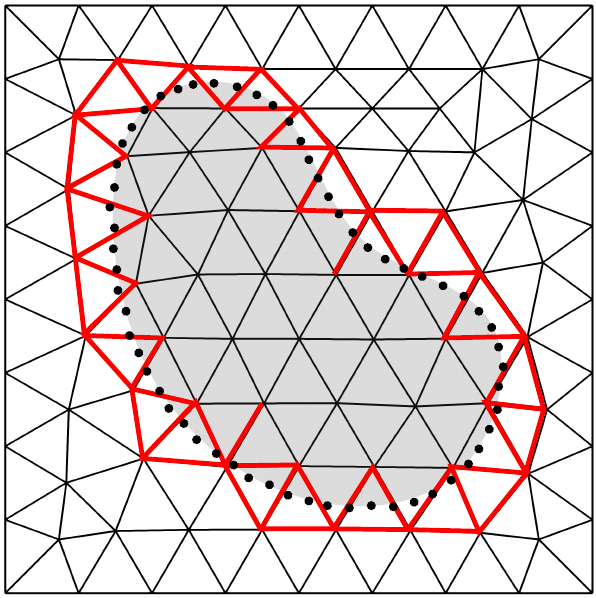}
    \caption{}
    \label{fig:imga_ghostPenaltyFaces}
  \end{subfigure}
  \caption{(a) An illustration of flow over an object represented by a point cloud. The object with boundary $\Gamma$ is immersed into the domain $\Omega$, which separates into a physical part $\Omega_{\text{phys}}$ and a fictitious part $\Omega_{\text{fict}}$. (b) Background mesh discretizing the domain $\Omega$. (c) The faces ($\mathcal{F}^f$) where the ghost penalty is applied are highlighted in red. These faces are shared by two elements that are fully or partially within the physical domain, where at least one of the elements is intersected by the boundary.}
\end{figure}

Let $\Omega \subset \mathbb{R}^d, \, d\in\left\{2,3\right\}$, be the computational domain, partitioned into two disjoint regions: the physical (fluid) domain $\Omega_{\text{phys}}$ and the fictitious domain $\Omega_{\text{fict}}$ separated by the immersed boundary $\Gamma$, as illustrated in Figure~\ref{fig:IMGA_domains}. The computational domain is discretized into a set of disjoint elements ${\Omega^e}$ such that $\Omega\subset\cup_e\overline{\Omega^e}$, as shown in Figure~\ref{fig:IMGA_mesh}. Since the mesh is not fitted to $\Gamma$, an element can be cut by the immersed boundary. We denote by $\Omega^e_\text{phys} = \Omega^e \cap \Omega_\text{phys}$ and $\Omega^e_\text{fict} = \Omega^e \cap \Omega_\text{fict}$ the portion of ${\Omega^e}$  belonging to the physical and fictitious domains respectively. The immersed boundary $\Gamma$, represented by a point cloud $\mathcal{P}$, is discretized into a collection of boundary elements ${\Gamma^b}$. Let $\mathcal{S}^h_u$ and $\mathcal{S}^h_p$ be the discrete trial function spaces for fluid velocity and pressure, respectively, and $\mathcal{V}^h_u$ and $\mathcal{V}^h_p$ be the corresponding test function spaces supported on ${\Omega^e}$. The IMGA discretization of the incompressible Navier--Stokes equations is formulated as follows: Find $\mathbf{u}^h \in \mathcal{S}^h_u$ and $p^h \in \mathcal{S}^h_p$ such that for all $\mathbf{w}^h \in \mathcal{V}^h_u$ and $q^h \in \mathcal{V}^h_p$,
\begin{align}\label{eq:weak_form}
\nonumber & B^\text{VMS}\left(\{\mathbf{w}^h, q^h\},\{\mathbf{u}^h,p^h\}\right) + B^\text{GhP}\left(\{\mathbf{w}^h, q^h\},\{\mathbf{u}^h,p^h\}\right) \\
& \quad\quad + B^\text{WBC}\left(\{\mathbf{w}^h, q^h\},\{\mathbf{u}^h,p^h\}\right)-F^\text{VMS}\left(\{\mathbf{w}^h, q^h\}\right) = 0\text{ ,} 
\end{align}
where the semi-linear form $B^\text{VMS}$ and the load vector $F^\text{VMS}$ are the terms associated with VMS, $B^\text{GhP}$ represents the ghost penalty stabilization, and $B^\text{WBC}$ includes the terms imposing the Dirichlet boundary conditions weakly. In this work, $B^\text{VMS}$ and $F^\text{VMS}$ are given as
\begin{align} \label{eq:B}
\nonumber & B^{\mathrm{VMS}}\left(\{\mathbf{w}^h, q^h\},\{\mathbf{u}^h, p^h\}\right)=\int_{\Omega_\text{phys}} \mathbf{w}^h \cdot \rho\left(\frac{\partial \mathbf{u}^h}{\partial t}+\mathbf{u}^h \cdot \boldsymbol{\nabla} \mathbf{u}^h\right) d \Omega +\int_{\Omega_\text{phys}} \boldsymbol{\varepsilon}(\mathbf{w}^h): \boldsymbol{\sigma}(\mathbf{u}^h, p^h) d \Omega \\
\nonumber & \quad +\int_{\Omega_\text{phys}} q^h \boldsymbol{\nabla} \cdot \mathbf{u}^h d \Omega -\sum_e \int_{\Omega^e_\text{phys}}\left(\mathbf{u}^h \cdot \boldsymbol{\nabla} \mathbf{w}^h+\frac{\boldsymbol{\nabla} q^h}{\rho}\right) \cdot \mathbf{u}^{\prime} d \Omega -\sum_e \int_{\Omega^e_\text{phys}} p^{\prime} \boldsymbol{\nabla} \cdot \mathbf{w}^h d \Omega \\
\nonumber & \quad+\sum_e \int_{\Omega^e_\text{phys}} \mathbf{w}^h \cdot\left(\mathbf{u}^{\prime} \cdot \boldsymbol{\nabla} \mathbf{u}^h\right) d \Omega -\sum_e \int_{\Omega^e_\text{phys}} \frac{\boldsymbol{\nabla} \mathbf{w}^h}{\rho}:\left(\mathbf{u}^{\prime} \otimes \mathbf{u}^{\prime}\right) d \Omega \\
& \quad+\sum_e \int_{\Omega^e_\text{phys}}\left(\mathbf{u}^{\prime} \cdot \boldsymbol{\nabla} \mathbf{w}^h\right) \bar{\tau} \cdot\left(\mathbf{u}^{\prime} \cdot \boldsymbol{\nabla} \mathbf{u}^h\right) d \Omega\text{ ,}
\end{align}
and
\begin{align}\label{eq:F}
F^\text{VMS}\left(\{\mathbf{w}^h, q^h\}\right) =&\int_{\Omega_\text{phys}}\mathbf{w}^h\cdot\rho\mathbf{f}~d\Omega 
+ \int_{\Gamma^{\text{N}}}\mathbf{w}^h\cdot\mathbf{h}~d\Gamma\text{ ,}
\end{align}
where $\rho$ is the fluid density, $\mathbf{f}$ is the body force per unit mass, $\mathbf{h}$ is the traction vector applied at the Neumann boundary $\Gamma^\text{N} \subset \Gamma$, $\pmb{\sigma}$ is the Cauchy stress tensor, and $\pmb{\varepsilon}$ is the strain-rate tensor. $\mbox{$\pmb{\sigma} (\mathbf{u}, p) = -p \, \mathbf{I} + 2\mu\,\pmb{\varepsilon}(\mathbf{u})$}$ and $\mbox{$\pmb{\varepsilon}(\mathbf{u}) = \frac{1}{2} (\nabla \mathbf{u}+ \nabla \mathbf{u}^T)$}$, where $\mathbf{I}$ is the identity tensor and $\mu$ is the dynamic viscosity. Additionally, the fine-scale velocity $\mathbf{u}^{\prime}$ and pressure $p^{\prime}$ are defined as
\begin{align}\label{eq:u_prime}
&\mathbf{u}^{\prime}=-\tau_{\mathrm{M}}\left(\rho\left(\frac{\partial \mathbf{u}^h}{\partial t}+\mathbf{u}^h \cdot \nabla \mathbf{u}^h-\mathbf{f}\right)-\boldsymbol{\nabla} \cdot \boldsymbol{\sigma}\left(\mathbf{u}^h, p^h\right)\right)\text{ ,} \\
\label{eq:p_prime}
&p^{\prime}=-\rho \tau_{\mathrm{C}} \boldsymbol{\nabla} \cdot \mathbf{u}^h\text{ ,}
\end{align}
where $\bar{\tau}$, $\tau_\text{M}$, and $\tau_\text{C}$ are the stabilization parameters, and their detailed expressions can be found in Ref.~\cite{Xu16tetra}. Other options for the stabilization parameters can be found in Refs.~\cite{Tezduyar00Finit, Hsu10Impro, Takizawa18Stabi, Jia23time, Takizawa23Varia}.

\subsection{Ghost penalty stabilization on cut elements}
\label{subsec:ghostPenalty}
In immersed methods, the boundary intersects the background mesh arbitrarily. This can create small-cut/sliver elements, where only a small fraction of the element lies in the physical domain. The basis functions associated with such small-cut elements have limited support within the physical domain and can lead to an ill-conditioned matrix system. The ghost penalty method~\cite{Burman10Ghost, Burman12Ficti, Burman15CutFE} adds an additional stabilization to the variational formulation near the cut elements and improves the conditioning of the matrix system.
This method has been applied across a wide range of problems and demonstrated improved stability and accuracy~\cite{Gurkan19stabi, Stoter23Criti, Divi22Resid, Massing14stabi, Schott14new, Schott15face, Schott16stabi, Liu21Nitsc}. In this work, we employ ghost penalty stabilization for both velocity and pressure in incompressible flow based on the work of Refs.~\cite{Schott14new, Dettmer16stabi}. The ghost penalty terms, represented by $B^\text{GhP}$ in Eq.~\eqref{eq:weak_form}, are given as follows:
\begin{align}\label{eq:ghostPenalty}
\nonumber B^\text{GhP}\left(\{\mathbf{w}^h, q^h\},\{\mathbf{u}^h,p^h\}\right) &= \sum_{f} \int_{\mathcal{F}^f} \tau^\text{GhP}_\mathbf{u} \llbracket \nabla_\mathbf{n} \mathbf{w}^h \rrbracket\cdot \llbracket \nabla_\mathbf{n} \mathbf{u}^h \rrbracket d \Gamma \\
&+ \sum_{f} \int_{\mathcal{F}^f} \tau^\text{GhP}_p \llbracket \nabla_\mathbf{n} q^h \rrbracket \llbracket \nabla_\mathbf{n} p^h \rrbracket d \Gamma\text{ ,}
\end{align}
where $\mathcal{F}^f$ are the faces shared by two elements within the physical domain $\Omega_\text{phys}$, in which at least one of the elements is cut, as illustrated in Figure~\ref{fig:imga_ghostPenaltyFaces}. $\llbracket \cdot \rrbracket$ denotes the jump operator evaluated at the face $\mathcal{F}^f$, and $\tau^\text{GhP}_\mathbf{u}$ and $\tau^\text{GhP}_p$ are the ghost penalty stabilization parameters \cite{Schott14new} given by 
\begin{align}
&\tau^\text{GhP}_\mathbf{u} = \alpha^\text{GhP}_\text{visc}\, \mu h_f + \alpha^\text{GhP}_\text{conv}\, \rho \, \|\mathbf{u}^h \cdot \mathbf{n} \| \, h_f^2 \text{ ,}\\ 
&\tau^\text{GhP}_p = \alpha^\text{GhP}_p\, \left( \frac{\mu}{h_f} + \rho \right)^{-1} h_f^2 \text{ .}
\end{align}
In the above, $h_f$ is the maximum length in the normal direction of each neighboring element at face $\mathcal{F}^f$, and $\alpha^\text{GhP}_\text{visc}$, $\alpha^\text{GhP}_\text{conv}$ and  $\alpha^\text{GhP}_p$ are non-dimensional penalty parameters. The choice of penalty parameter is a trade-off between the stability and accuracy of the solution. $\alpha^\text{GhP}_\text{visc}$, $\alpha^\text{GhP}_\text{conv}$, $\alpha^\text{GhP}_p \in [0.001, 0.05]$ are recommended for incompressible flow~\cite{Schott14new, Schott15face, Schott16stabi} and we found $\alpha^\text{GhP}_\text{visc}$, $\alpha^\text{GhP}_\text{conv}$, $\alpha^\text{GhP}_p = 0.05$ effective for all cases considered in this study. In Eq.~\eqref{eq:ghostPenalty}, the ghost penalty method penalizes the jump of the normal gradient of the solution within the cut elements. This enforces smoothness and provides stability to both the velocity and pressure fields across the cut elements.
 
\subsection{Non-symmetric Nitsche-based weak BC}
\label{sec:weakbc}
In contrast to boundary-fitted methods, where the Dirichlet no-slip condition is applied strongly on the nodes, IMGA imposes the no-slip condition directly on the immersed surface via Nitsche-based weak BC~\cite{Bazilevs07Weak1, Ruberg12Subdi, Kamensky15immer}. A comprehensive review of weakly enforced Dirichlet BCs in computational flow analysis, including their role in IMGA, can be found in Hsu et al.~\cite{Hsu26revie}. The weak BC terms, represented by $B^\text{WBC}$ in Eq.~(\ref{eq:weak_form}), are given as
\begin{align}\label{eq:incomp_wbc}
\nonumber & B^\text{WBC}\left(\{\mathbf{w}^h, q^h\},\{\mathbf{u}^h,p^h\}\right) = 
-\sum_b\int_{\Gamma^b \bigcap\Gamma^\text{D}} \mathbf{w}^h \cdot\left(-p^h \mathbf{n}+2 \mu\,\pmb{\varepsilon}(\mathbf{u}^h)\, \mathbf{n}\right) d \Gamma \\
\nonumber & \quad -\sum_b\int_{\Gamma^b \bigcap\Gamma^\text{D}}q^h \mathbf{n}\cdot\left(\mathbf{u}^h-\mathbf{g}\right) d \Gamma
-\sum_b\int_{\Gamma^b \bigcap\Gamma^\text{D}} \gamma \,\left(2 \mu\,\pmb{\varepsilon}(\mathbf{w}^h)\, \mathbf{n}\right) \cdot\left(\mathbf{u}^h-\mathbf{g}\right) d \Gamma \\
\nonumber & \quad-\sum_b\int_{\Gamma^b \bigcap\Gamma^{\text{D},-}} \mathbf{w}^h \cdot \rho\left(\mathbf{u}^h \cdot \mathbf{n}\right)\left(\mathbf{u}^h-\mathbf{g}\right) d \Gamma \\
\nonumber & \quad +\sum_b\int_{\Gamma^b \bigcap\Gamma^\text{D}} \tau_{\mathrm{TAN}}^B\left(\mathbf{w}^h-\left(\mathbf{w}^h \cdot \mathbf{n}\right) \mathbf{n}\right) \cdot\left(\left(\mathbf{u}^h-\mathbf{g}\right)-\left(\left(\mathbf{u}^h-\mathbf{g}\right) \cdot \mathbf{n}\right) \mathbf{n}\right) d \Gamma \\
 & \quad+\sum_b\int_{\Gamma^b \bigcap\Gamma^\text{D}} \tau_{\mathrm{NOR}}^B\left(\mathbf{w}^h \cdot \mathbf{n}\right)\left(\left(\mathbf{u}^h-\mathbf{g}\right) \cdot \mathbf{n}\right) d \Gamma\text{ ,}
\end{align}
where $\mathbf{n}$ is the unit normal, $\Gamma^{\text{D},-}$ is the inflow part of the Dirichlet boundary $\Gamma^\text{D}$, on which $\mathbf{u}^h\cdot\mathbf{n} < 0$, and $\tau_{\mathrm{TAN}}^B$ and $\tau_{\mathrm{NOR}}^B$ are stabilization parameters that act on the tangential and normal components of the velocity at the boundary, respectively. The first term is the consistency term obtained from integration by parts of the Cauchy stress of the Navier--Stokes equations. The next two terms are the adjoint consistency terms corresponding to the pressure and the viscous term, where $\gamma=1$ and $\gamma = -1 $ yield the symmetric and non-symmetric variants of Nitsche's method, respectively. The fourth term stabilizes the convective instability on portions of the immersed boundary where the flow enters the domain~\cite{Bazilevs07Weak1}. The fifth and sixth terms are the stabilization terms in the tangential and normal directions, respectively. Previous IMGA and point-cloud CFD studies employed the symmetric variant ($\gamma= 1$) with sufficiently large constant stabilization parameters, which produce accurate fields and integrated forces~\cite{Kamensky15immer, Xu16tetra, Balu23Direc, Hsu16Direc, Wang23Photo, Jaiswal24Mesh, Corpuz25Direc}. However, as discussed in Section~\ref{sec:intro}, this choice is inadequate for computing the local WSS distribution. A constant large parameter contributes a non-smooth, amplified residual to the boundary traction. Moreover, estimating the parameters through an element-local eigenvalue problem~\cite{Embar10Impos} to compute the lower bound necessary for coercivity produces cut-sensitive, spurious values that affect the local traction distribution. This behavior is demonstrated in Section~\ref{sec:validation}. Therefore, we employ the non-symmetric variant ($\gamma = -1$), which removes the strict lower bound on the stabilization parameters, while retaining the stabilization terms for accuracy and robustness necessary for complex problems~\cite{Guo17param, Heimann13unfit}. In this study, the parameters are chosen based on physical characteristics of the flow, independent of the cut configuration.

\subsubsection{Tangential stabilization parameter}
\label{sec:weakbc_tangent}
The tangential stabilization term allows the flow to slip on the wall, with the amount of slip controlled by $\tau^B_\text{TAN}$, where large values approach the strong no-slip and smaller values admit more slip. This is closely related to near-wall modeling techniques in turbulence modeling~\cite{Pope00Turbu}, where the strong no-slip condition is replaced by a wall-stress condition determined with wall models. The connection between weakly enforced no-slip conditions and traditional wall models was established by Bazilevs et al.~\cite{Bazilevs07Weak2} and further developed by Golshan et al.~\cite{Golshan15Large}, where wall-model-based parameters showed improved performance on under-resolved meshes compared to the standard mesh-size-based scaling. Following this approach, we employ a near-wall model to estimate the tangential stabilization parameter. We note that in Ref.~\cite{Golshan15Large}, the wall-modeled parameter had to be bounded below by the value from the inverse estimate to satisfy the coercivity requirement for the symmetric Nitsche formulation. The non-symmetric variant imposes no such constraint, so the parameter is determined solely by the wall model, providing a controlled flow slip that reflects the near-wall physics. In near-wall modeling, the WSS is expressed as
\begin{align} \label{eq:trad_wss}
\boldsymbol{\tau}_\text{wm} = \rho u^{*2} \frac{\mathbf{u}_t - \mathbf{g}_t}{\|\mathbf{u}_t- \mathbf{g}_t\|} \text{ ,}
\end{align}
where $u^{*}$ is the wall-friction velocity magnitude and $\mathbf{u}_t = \mathbf{u}^h - (\mathbf{u}^h\cdot \mathbf{n})\,\mathbf{n} $ and $\mathbf{g}_t = \mathbf{g} - (\mathbf{g}\cdot \mathbf{n})\,\mathbf{n} $. The wall-friction velocity $u^{*}$ is computed iteratively through Spalding's single formula parameterization of the turbulent layer~\cite{Spalding61singl}
\begin{align}\label{eq:spalding}
y^{+} = u^{+} + e^{-\kappa B}\left( e^{\kappa u^{+}} -1 -\kappa u^{+} -\frac{(\kappa u^{+})^2}{2} -\frac{(\kappa u^{+})^3}{6}   \right) \text{ ,}
\end{align}
where
\begin{align}
\nonumber y^{+} = \frac{y^{\text{A}} u^* \rho}{\mu} \text{ ,} \quad u^{+} = \frac{\|\mathbf{u}^{\text{A}}_t - \mathbf{g}_t\|}{u^*} \text{ ,}\quad \kappa = 0.41 \text{,} \text{ and } B = 5.5 \text{ .}
\end{align}
Here, $y^\text{A}$ is the wall-normal distance to a sampling point A, and $\mathbf{u}^{\text{A}}_t$ is the tangential velocity sampled there. In traditional near-wall modeling, the sampling point A is placed a few grid points off the wall, such that it lies within the log region of the turbulent boundary layer~\cite{Kawai12Wallm}. In this work, we place point A at a fixed wall-normal distance equal to the average wall-normal size of the cut elements. This construction allows for consistent $y^\text{A}$ across all discrete boundary points, and it positions the sampling location one element beyond the cut element, so that it samples the velocity from elements unaffected by the cut. 

The wall-model shear stress in Eq.~\eqref{eq:trad_wss} is precisely the tangential traction that the weak BC should exert on the wall. The viscous term contributes the resolved shear stress, while the stabilization term provides the unresolved shear stress, analogous to a subgrid-scale model. This relation can be expressed as
\begin{align}
\label{eq:wm_bal}
\boldsymbol{\tau}_\text{wm} = - 2 \mu\,\pmb{\varepsilon}(\mathbf{u}^h)\, \mathbf{n} + \tau_{\mathrm{TAN}}^B \left(\mathbf{u}_t - \mathbf{g}_t\right) \text{ .}
\end{align}
Approximating the viscous traction by its dominant wall-normal component, equating Eq.~\eqref{eq:trad_wss} and Eq.~\eqref{eq:wm_bal}, and solving for the parameter gives
\begin{align}\label{eq:tau_tan}
\tau_{\mathrm{TAN}}^B = \frac{1}{\|\mathbf{u}_t - \mathbf{g}_t\|}\left(\rho\,u^{*2} - \mu \left( \frac{\partial\|\mathbf{u}_t\|}{\partial y} \right)_{\Gamma^b}  \right) \text{ ,}
\end{align}
where
\begin{align}\label{eq:dudn}
\left( \frac{\partial\|\mathbf{u}_t\|}{\partial y} \right)_{\Gamma^b} \approx \frac{\|\mathbf{u}^{\text{A}}_t\| - \|\mathbf{u}_t\|}{y^{\text{A}}} \text{ .}
\end{align}

\subsubsection{Normal stabilization parameter}
\label{sec:weakbc_normal}
The normal stabilization parameter controls the enforcement of the no-penetration boundary condition at the wall. In contrast to the tangential direction, where the stabilization allows the flow to slip slightly, the normal condition is a kinematic constraint that must be enforced strictly. Any normal velocity residual corresponds to spurious mass flux through the immersed boundary. Motivated by the scaling of stabilization parameters in the VMS formulation, we employ
\begin{align}
\tau_{\mathrm{NOR}}^B= C_b \left( \frac{\mu}{h_n} + \rho \|\mathbf{u}^h\| + \frac{\rho h_n}{\Delta t} \right) \text{ ,}
\end{align}
where $h_n$ is the wall-normal element size, and $C_b$ is a non-dimensional constant. In this work, $C_b = 10^2$ is found to be effective for all cases considered in this study. The three terms scale with the viscous, convective, and transient character of the flow, so that the stabilization adjusts with dominant physics across flow regimes.  

\subsection{Wall shear stress recovery}
\label{sec:wss_recovery}
In the finite element method, stresses are computed from the gradient of the discrete solution and are naturally evaluated at the quadrature points within each element. The resulting stress field is discontinuous across element boundaries, and piecewise constant for linear elements. For post-processing purposes, including visualization, a posteriori error estimation, and adaptive mesh refinement~\cite{Zienkiewicz87simpl}, the stress is commonly projected to a continuous field at the nodes (or control points, as in Isogeometric Analysis (IGA)~\cite{Hughes05Isoge}), using techniques such as nodal averaging, lumped-mass projection, or global $L^2$ projection~\cite{Hinton74Local}. Among these, patch-based recovery methods have become popular owing to their superconvergence properties and are implemented in commercial finite element software such as Abaqus and ANSYS~\cite{Guo25Recov}. The superconvergent patch recovery (SPR) by Zienkiewicz and Zhu~\cite{Zienkiewicz92super} constructs a patch of elements around each node and fits a polynomial, in the least-squares sense, to the stresses sampled at the superconvergent points of each element, yielding recovered nodal stresses of higher accuracy than the directly computed field. Element-based variants, proposed by Wiberg and Li~\cite{Wiberg94Super}, form a patch around an element instead, and obtain nodal stresses by weighted combination of polynomials from overlapping patches. Subsequent development enhanced the recovery by constraining it with equilibrium and boundary conditions~\cite{Blacker94Super}. These methods have a particular advantage in immersed discretizations, where small cut elements are known to produce spurious, non-physical stress concentrations. By gathering a patch of neighboring elements, the recovered field regularizes the cut-induced artifacts that pointwise evaluation cannot avoid~\cite{Zhang22impro, Navarro20Super, Capatina21Flux}. In this work, we build on this foundation and propose a wall shear stress recovery method for point cloud-based CFD, in which the recovered stress field is constructed on a geometry-following element-based patch and additionally constrained by the variationally consistent traction from the weak BC. 

\subsubsection{Element-based patch construction}
\label{sec:elementpatch}
\begin{figure}[!t]
    \centering
    \includegraphics[width=\textwidth]{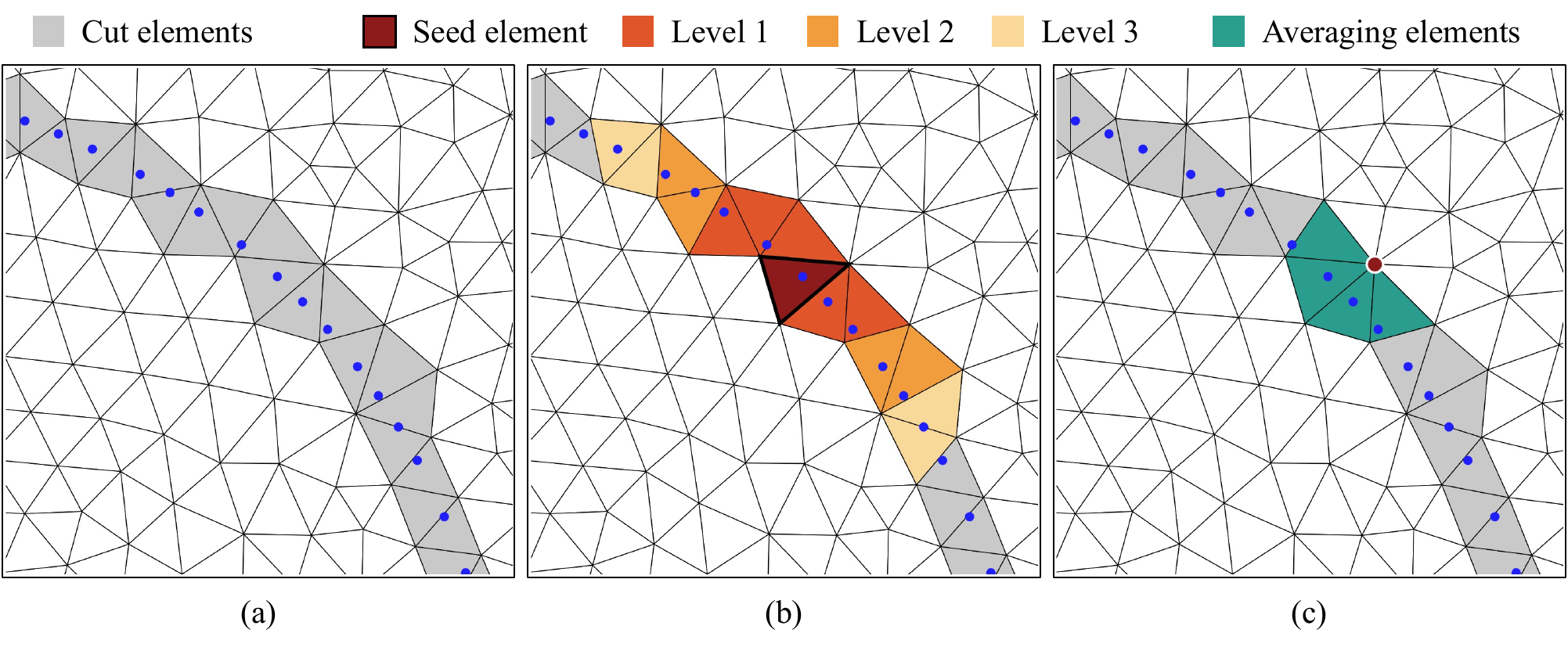}
    \caption{Geometry-following element-based patch construction for WSS recovery, illustrated on a two-dimensional background mesh. (a) Point cloud representing the boundary (blue) and the cut elements (gray). (b) Three-level patch of a seed element, restricted to node-connected cut elements. (c) Nodal averaging over the cut elements sharing the highlighted node, each contributing a recovered field from its own patch.}
    \label{fig:patch}
\end{figure}

We are interested in evaluating the recovered stress field of cut elements at the points of the immersed point cloud to obtain the WSS. To this end, we employ a geometry-following element-based patch: the patch of a cut element consists of its neighboring elements that are themselves cut by the geometry. Because the cut-element layer geometrically follows the immersed boundary, this approach confines the polynomial fit to the near-wall region and prevents steep boundary-layer gradients from being diffused by low-gradient fields farther from the wall. In contrast to classical SPR, where the polynomial is fitted to stresses sampled at the superconvergent points, we determine the recovered field by a least-squares fit in the integral sense over the physical portion of the patch elements. This choice is naturally suited to the immersed setting where each element contributes in proportion to its physical volume, such that small cut elements, whose directly computed stresses are the least reliable, carry correspondingly smaller weight in the recovery. The patch reconstruction is summarized in Algorithm~\ref{alg:patch}, where the extent of the patch is controlled by the connectivity level as shown in Figure~\ref{fig:patch}. The conditioning of the local least-squares system depends on the size and shape of the integration domain. On a boundary-fitted mesh, a one-level patch provides adequate support in the interior, while nodes on the boundary are classically treated using adjacent interior patches, whose support is otherwise insufficient~\cite{Zienkiewicz92super}. In our geometry-following patch, the cut elements form a thin, quasi-two-dimensional layer along the boundary, and their physical portion may occupy only a fraction of each element. We therefore use three levels, which were found to yield integration domains comparable in extent to one-level interior patches and a well-conditioned fit. A similar patch enlargement procedure has been employed in Navarro-Jim{\'e}nez et al.~\cite{Navarro20Super}, where patches with a low material-to-element volume ratio are extended to gather sufficient support.

\begin{algorithm}[!t]
    \caption{Geometry-following patch construction for a cut element $e_0$}
    \label{alg:patch}
    \begin{algorithmic}[1]
    \State $\mathrm{patch} \gets \{\, e_0 \, \}$
    \State $\mathrm{frontier} \gets \{\, e_0 \, \}$
    \For{$\ell = 1, \dots, \mathrm{level}$}
        \State $\mathrm{next} \gets \emptyset$
        \For{\textbf{each} $e \in \mathrm{frontier}$}
            \For{\textbf{each} elem $e'$ sharing a node with $e$}
                \If{$e'$ is cut \textbf{and} $e' \notin \text{patch}$}
                    \State $\mathrm{patch} \gets \mathrm{patch} \, \cup \, \{{e'\}} $
                    \State $\mathrm{next} \gets \mathrm{next} \, \cup \, \{{e'\}} $
                \EndIf
            \EndFor
        \EndFor
        \State $\mathrm{frontier} \gets \mathrm{next}$
    \EndFor
    \State \Return $\mathrm{patch}$
    \end{algorithmic}
\end{algorithm}

\subsubsection{Recovery}
For incompressible flow, the Cauchy stress tensor decomposes into volumetric (pressure) and deviatoric (viscous) contributions,
\begin{align}
\boldsymbol{\sigma}^h = -p^h\mathbf{I} + \mathbf{s}^h \text{ ,}
\end{align}
where $\mathbf{s}^h = 2 \mu\,\pmb{\varepsilon}(\mathbf{u}^h)$. We apply the recovery exclusively to the deviatoric part, $\mathbf{s}^h$, because the WSS is a purely tangential quantity. As a result, it depends entirely on the viscous term, whereas the pressure contribution acts solely in the wall-normal direction. Using the Voigt notation, the symmetric deviatoric stress tensor is represented as a column vector,
\begin{align}
\mathbf{s}^h = 
        \begin{bmatrix}
        s^h_{1} &
        s^h_{2} &
        s^h_{3} &
        s^h_{4} &
        s^h_{5} &
        s^h_{6} 
        \end{bmatrix}^T
        =
        \begin{bmatrix}
        s^h_{xx} &
        s^h_{yy} &
        s^h_{zz} &
        s^h_{yz} &
        s^h_{xz} &
        s^h_{xy} 
        \end{bmatrix}^T \text {.}
\end{align}
The stress recovery is performed locally over $\Omega^p$, which defines the physical domain of an element-based patch for a given cut element, as detailed in Section~\ref{sec:elementpatch}. Over each patch, we construct a recovered stress field $\mathbf{s}^*$ using a linear polynomial basis
\begin{align}
\mathbf{P}(\mathbf{x}) = \begin{bmatrix} 1 & x & y & z \end{bmatrix} \in \mathbb{R}^{1\times 4}  \text{ ,}
\end{align}
and each stress component is approximated as 
\begin{align}
    s^*_{i}(\mathbf{x}) = \mathbf{P}(\mathbf{x})\,\mathbf{a}_{i} \text{ ,} \quad i = 1,\ldots, 6,
\end{align}
where $\mathbf{a}_{i} \in \mathbb{R}^{4\times 1}$ is the column vector of unknown coefficients of the $i$-th stress component. Collecting all six components in block matrix form gives
\begin{align}
    \mathbf{s}^*(\mathbf{x}) = \mathbf{M}(\mathbf{x}) \,\mathbf{a} \text{ ,}
\end{align}
where 
\begin{align}
\mathbf{M} = \operatorname{diag}\left( \mathbf{P}, \mathbf{P}, \mathbf{P}, \mathbf{P}, \mathbf{P}, \mathbf{P} \right)  \in \mathbb{R}^{6\times24} \text{ ,}
\end{align}
and
\begin{align}
\mathbf{a}= \begin{bmatrix} \mathbf{a}_{1}^T & \mathbf{a}_{2}^T & \mathbf{a}_{3}^T & \mathbf{a}_{4}^T & \mathbf{a}_{5}^T & \mathbf{a}_{6}^T \end{bmatrix}^T \in \mathbb{R}^{24\times 1} \text{ .}
\end{align}
In standard patch recovery, the coefficients $\mathbf{a}$ are determined by minimizing the least-squares error between the recovered stress and finite element stress over the patch,
\begin{align}
    \label{eq:spr_internal_functional}
    J_\Omega(\mathbf{a}) = \frac{1}{2} \int_{\Omega^p} \left\| \mathbf{M}\mathbf{a} - \mathbf{s}^h \right\|^2 \, d\Omega \text{ .}
\end{align}
Although this minimization produces a smooth stress field, its accuracy degrades near the immersed boundary because it fails to enforce traction compatibility. Specifically, purely interior recovery neglects the additional stabilization terms that are incorporated alongside the standard Cauchy stress to define a variationally consistent traction with weak BC,
\begin{align}
    \label{eq:var_traction}
    \nonumber \mathbf{t}^h = &- \boldsymbol{\sigma}^h \mathbf{n}- \rho \left[ \mathbf{u}^h \cdot \mathbf{n}\right]_{-} \left(\mathbf{u}^h - \mathbf{g} \right) \\
    &+ \tau^B_\text{TAN} \left(\left(\mathbf{u}^h-\mathbf{g}\right)-\left(\left(\mathbf{u}^h-\mathbf{g}\right) \cdot \mathbf{n}\right) \mathbf{n}\right) + \tau^B_\text{NOR}\left(\left(\mathbf{u}^h-\mathbf{g}\right) \cdot \mathbf{n}\right) \mathbf{n} \text{ .}
\end{align}
The tangential component of this traction,
\begin{align}
    \mathbf{t}^h_t =  \mathbf{t}^h - \left(\mathbf{t}^h \cdot \mathbf{n} \right) \, \mathbf{n} \text{ ,}
\end{align}
is the finite element WSS, and it depends entirely on the deviatoric stress, which is consistent with restricting the recovery to $\mathbf{s}^h$. We therefore augment the functional~\eqref{eq:spr_internal_functional} with a boundary term that penalizes the difference between the traction of the recovered stress and $\mathbf{t}^h_t$ over the portion of the immersed boundary covered by the patch, $\Gamma^p = \Gamma^b \cap \overline{\Omega^p} $:
\begin{align}
    \label{eq:spr_full_functional}
    J(\mathbf{a}) = \frac{1}{2}\int_{\Omega^p} \left\| \mathbf{M}\mathbf{a} - \mathbf{s}^h \right\|^2 \, d\Omega + \frac{\alpha}{2}\int_{\Gamma^p} \left\| \mathbf{Q}\mathbf{M}\mathbf{a} - \mathbf{t}^h_t \right\|^2 \, d\Gamma \text{ ,}
\end{align}
where $\alpha$ is a penalty parameter, and $\mathbf{Q} \in \mathbb{R}^{3\times6}$ defines the Cauchy relation mapping the Voigt stress vector to the corresponding boundary traction,
\begin{align}
    \mathbf{Q} = -\begin{bmatrix}
        n_x & 0 & 0 & 0 & n_z & n_y \\
        0 & n_y & 0 & n_z & 0 & n_x \\
        0 & 0 & n_z & n_y & n_x & 0
    \end{bmatrix} \text{ ,}
\end{align}
where $n_x$, $n_y$, $n_z$ are components of the outward unit normal $\mathbf{n}$ on the boundary. The boundary term follows the classical enhancements of SPR, in which weighted residuals of the traction boundary conditions are added to the least-squares functional~\cite{Blacker94Super, Wiberg95Impro}. Later improvements, such as SPR-C~\cite{Rodenas07Impro}, impose the boundary traction as a hard constraint through Lagrange multipliers, which corresponds to the limit $\alpha \rightarrow \infty$. In this work, we retain the penalty form and set $\alpha = 1$. Notably, under mesh refinement, the boundary integral scales as $\mathcal{O}(1/h)$ relative to the interior domain integral. As a result, the boundary traction constraint is enforced progressively more strongly, a desirable property given that the weak BC traction is the most reliable boundary traction in this immersed approach. The minimizer of the modified functional (Eq.~\eqref{eq:spr_full_functional}) is obtained by setting $\partial J /\partial \mathbf{a} = \mathbf{0}$ leading to the linear system,
\begin{align}
    \left( \mathbf{A}_\Omega + \alpha \mathbf{A}_\Gamma \right) \, \mathbf{a} = \mathbf{b}_\Omega + \alpha \mathbf{b}_\Gamma \text{ ,}
\end{align}
where 
\begin{align}
    \nonumber \mathbf{A}_\Omega &= \int_{\Omega^p} \mathbf{M}^T\mathbf{M} \, d\Omega \text{ ,} \quad 
    \mathbf{A}_\Gamma = \int_{\Gamma^p} \left( \mathbf{Q}\mathbf{M} \right)^T\left( \mathbf{Q}\mathbf{M} \right) \, d\Gamma \text{ ,} \\
    \nonumber \mathbf{b}_\Omega &= \int_{\Omega^p} \mathbf{M}^T\mathbf{s}^h \, d\Omega \text{ ,} \quad
    \mathbf{b}_\Gamma = \int_{\Gamma^p} \left( \mathbf{Q}\mathbf{M} \right)^T\mathbf{t}^h_t \, d\Gamma \text{ .}
\end{align}
Here, $ \mathbf{A}_\Omega$ is block-diagonal, coupling the six stress components only through the boundary term. This generates a $24 \times 24$ symmetric positive-definite linear system per patch, which is solved directly via LU factorization. Once the coefficients $\mathbf{a}$ are determined for each cut element, the recovered polynomial is evaluated at the element nodes. As a given node is generally shared by multiple cut elements, it accumulates distinct recovered fields from each adjacent element. These contributions are subsequently averaged (Figure~\ref{fig:patch}c), consistent with standard element-based patch recovery procedures~\cite{Wiberg95Impro, Zhang22impro}. This process produces a continuous recovered deviatoric stress field defined over the background mesh. Finally, the recovered wall shear stress, $\boldsymbol{\tau}_w$, at each discrete boundary point is computed by interpolating this continuous field at the point's location within its corresponding background element, applying the Cauchy stress relation, and extracting the tangential component:
\begin{align}
    \boldsymbol{\tau}_w = \mathbf{t}^* - \left( \mathbf{t}^* \cdot \mathbf{n} \right) \, \mathbf{n} \text{ ,} \quad \mathbf{t}^* = \mathbf{Q}\, \mathbf{s}^*(\mathbf{x}) \text{ ,}
\end{align}
where $\mathbf{s}^*(\mathbf{x})$ denotes the interpolated recovered stress at the point $\mathbf{x}$ and $\mathbf{n}$ is the outward unit normal at that point.

\subsection{Cut quadrature}
\label{sec:cut-quad}
\begin{figure}[!t]
    \centering
    \includegraphics[width=\textwidth]{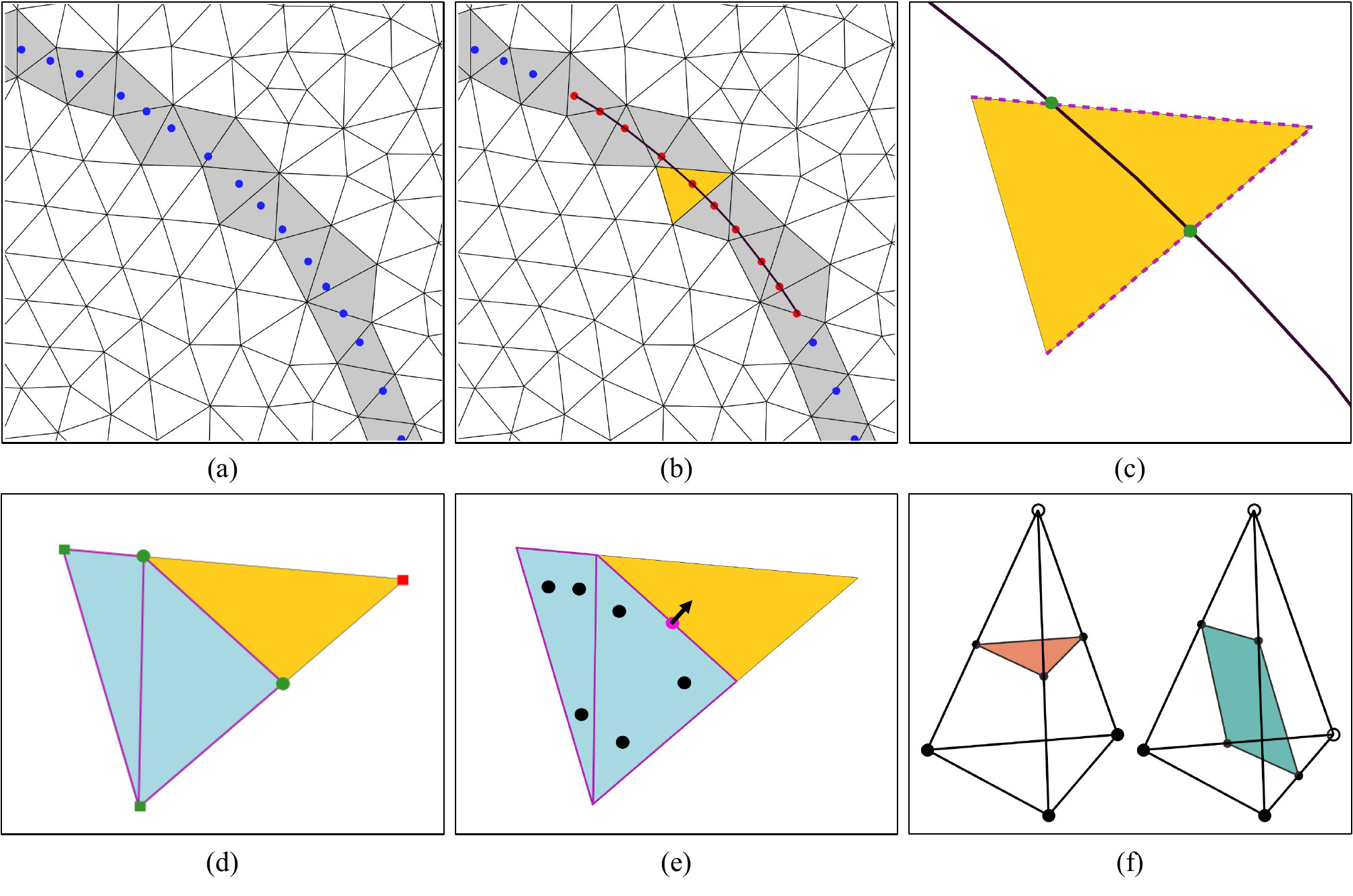}
    \caption{(a) Cut elements (gray) identified with winding number test on the background mesh. (b) Selection of a cut element and determining the local surface fit to its nearest points. (c) Intersecting the crossing edges with the local surface fit. (d) Intersected points and active vertices determine a sub-cell decomposition. (e) Standard quadrature rules are applied on the sub-cell, and a point with normal and area is sampled from the intersecting facet. (f) Illustration of possible intersection cases in a tetrahedral element: a triangle separating one or three vertices (left) and a quadrilateral separating two vertices (right). }
    \label{fig:cut-quad}
\end{figure}

Accurate integration over the intersected elements is essential for the accuracy of the IMGA framework: the VMS formulation is integrated over the physical domain $\Omega_\text{phys}$, the weak BC over the immersed boundary $\Gamma^b$, and the stress recovery over both. Standard quadrature rules are constructed for smooth integrands over full elements and lose accuracy when the integrand is restricted to an arbitrarily shaped physical portion. To address this, IMGA and the finite cell method employ a subdivision-based algorithm that adaptively refines the quadrature resolution near the boundary to construct the volume integration rule~\cite{Kamensky15immer, Xu16tetra, Varduhn16tetra, Balu23Direc, Duczek16finit, Kudela13Highl}. This strategy is simple and robust, but the number of quadrature points grows rapidly with subdivision level. Furthermore, accurate surface integration requires that every intersected element contains at least one integration point~\cite{Hsu16Direc, Balu23Direc}. To achieve this in point cloud-based CFD, a mesh-driven resampling method~\cite{Jaiswal24Mesh} was proposed for real-world scans, which come in a fixed number and density, to resample the points such that every intersected background element has one point for surface integration. The method relied on winding number-based inside/outside classification~\cite{Jacobson13Robus, Barill18Fast} and intersecting the element with a local surface fit~\cite{Cazals05Estim} to resample a point. In this work, we propose a simple extension to this approach where the local fit is utilized to recover the cut portion of the element, as illustrated in Figure~\ref{fig:cut-quad}. 

First, a winding number test is performed to classify background nodes as inside or outside and then identify the intersected elements (Figure~\ref{fig:cut-quad}a). For each intersected element, a local quadratic surface is fitted using its $k$-nearest neighbors~\cite{Corpuz25Direc} in the point cloud (Figure~\ref{fig:cut-quad}b). Unlike the original mesh-driven resampling method, which intersected the element's median lines with the local surface to obtain a resampled point, we compute the intersection of the background element's crossing edges with the local fit (Figure~\ref{fig:cut-quad}c). The classified background nodes and the resulting edge intersections then determine a sub-cell decomposition of the element by a template-based sub-tetrahedralization (Figure~\ref{fig:cut-quad}d). Assuming a monotone boundary cut through a tetrahedral element, only a small number of intersection topologies exist: the cut separates either one or three vertices forming a triangular intersection or two vertices forming a quadrilateral (Figure~\ref{fig:cut-quad}f). Therefore, sub-cells are obtained by mapping a precomputed template rather than by general-purpose tessellation. Such marching-tetrahedra-type sub-decompositions are well established in immersed finite element frameworks~\cite{Burman15CutFE, Rim25Enric1, Rim25Enric2}. Finally, standard quadrature rules applied to the resulting sub-tetrahedra provide the volume quadrature for $\Omega_\text{phys}$, while sampling the intersection polygon with points, normals, and area provides the surface quadrature for $\Gamma^b$ within the element (Figure~\ref{fig:cut-quad}e). We emphasize that the local fits are constructed element-wise from unorganized points and used strictly for edge intersection without any global surface reconstruction, so the framework continues to operate directly on the point cloud. Moreover, the algorithm is more efficient than subdivision-based quadrature, as the inside/outside classification is performed once per background node rather than per quadrature point, and the sub-cell construction replaces the recursive subdivision entirely. 

\section{Validation}
\label{sec:validation}
This section validates the proposed method for performing direct flow simulation on point clouds and extracting accurate flow quantities, including WSS, using two canonical benchmarks: internal Hagen--Poiseuille pipe flow and external flow past a sphere. 

\subsection{Hagen--Poiseuille pipe flow}
\label{sec:hagen}
In this study, we simulate laminar flow within a cylindrical pipe, which is governed by the Hagen--Poiseuille law and provides an exact analytical solution for the axial pressure drop and WSS. For this problem,  we consider a pipe of length $L=5.0$ and radius $R=0.5$. A parabolic velocity profile, $u(r) = U_\text{max}(1-r^2/R^2)$, is prescribed at the inlet with centerline velocity $U_\text{max}= 1.0$, where $r$ is the radial distance from the centerline. A traction-free condition is applied at the outlet. The density and dynamic viscosity are set to $\rho = 1.0$ and $\mu = 0.01$, respectively, yielding a Reynolds number of $Re = \rho \bar{U} D/ \mu = 50$ based on the pipe diameter $D=2R$ and the mean velocity $\bar{U} = U_\text{max}/2$. For fully developed Hagen--Poiseuille flow, the pressure drop over the pipe length and the wall shear stress are given by
\begin{align}
    \nonumber \Delta p = \frac{4\mu L U_\text{max}}{R^2} = 0.8 \text{ ,} \quad \tau_w = \frac{2\mu U_\text{max}}{R} = 0.04 \text{ ,}
\end{align}
which serve as the exact reference values for the validation study.

\begin{figure}[!t]
    \centering
    \includegraphics[width=\textwidth]{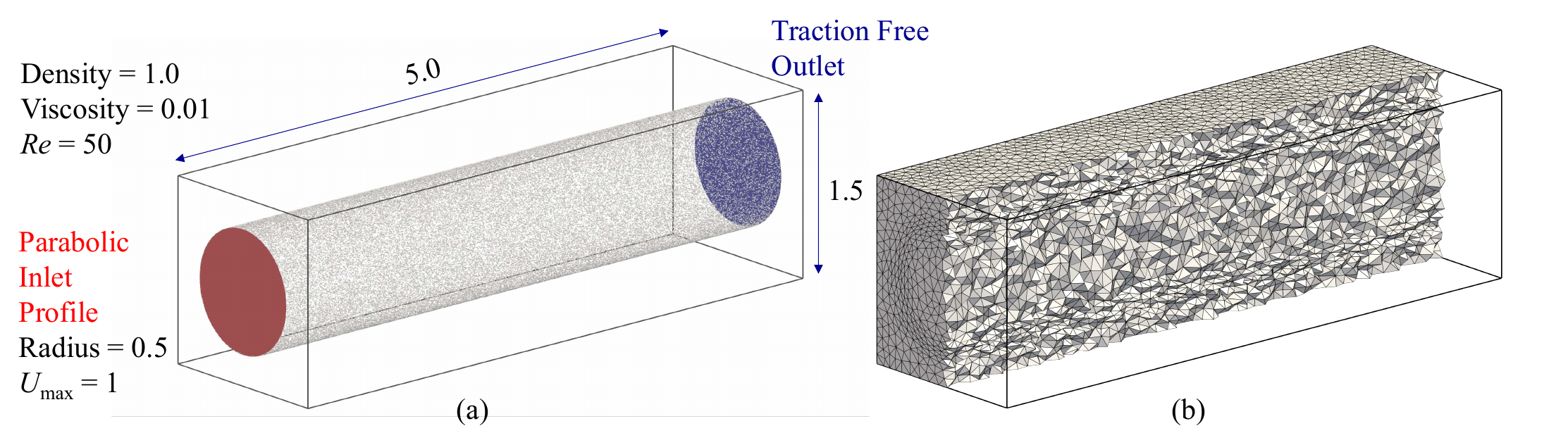}
    \caption{(a) Computational domain and (b) central cross-section of the background mesh IM0 considered for Hagen–Poiseuille pipe flow.}
    \label{fig:pipe-setup}
\end{figure}

\begin{table}[!t]\centering\small
    \centering
    \captionsetup{justification=centering}
    \newcommand{\tabincell}[2]{\begin{tabular}{@{}#1@{}}#2\end{tabular}}
    \caption{Background mesh statistics for the Hagen–Poiseuille pipe flow.}
    \begin{tabular}{lccccr}
        \toprule
        \multirow{2}{*}{Mesh} 
            & \multicolumn{3}{c}{Near wall} 
            & \multirow{2}{*}{Domain} 
            & \multirow{2}{*}{\tabincell{c}{Number of \\ active elements}} \\
        \cmidrule(lr){2-4}
            & $h_\mathrm{NOR}$
            & $h_\mathrm{TAN}$
            & \tabincell{c}{BL height} 
            & & \\
        \midrule
        IM0 & $0.04$  & $0.1$   & $0.04$ & $0.1$   & $80{,}219$    \\
        IM1 & $0.02$  & $0.05$  & $0.04$ & $0.05$  & $525{,}185$  \\
        IM2 & $0.01$  & $0.025$ & $0.04$ & $0.025$ & $3{,}362{,}170$ \\
        \bottomrule
    \end{tabular}
    \label{tab:meshstatpipe}
\end{table}

A point cloud representing a cylindrical surface of length $L=5.0$ and radius $R=0.5$ is generated by randomly sampling $10^5$ points from the analytical definition and is immersed into a computational domain, as shown in Figure~\ref{fig:pipe-setup}a. The domain fully encloses the point cloud laterally while conforming to the inlet and outlet planes, so that the parabolic inflow profile and traction-free BC can be applied consistently to replicate the exact solution. The no-slip condition is enforced weakly on the immersed boundary using the non-symmetric Nitsche-based weak BC ($\gamma = -1$), as detailed in Section~\ref{sec:weakbc}. The background domain is discretized into tetrahedral elements with refinement near the boundary, as presented in Figure~\ref{fig:pipe-setup}b. The domain is first meshed with Gmsh~\cite{Geuzaine09Gmsh} and subsequently re-meshed with anisotropic refinement along the boundary using Mmg3D~\cite{Dapogny14Three}, following the procedure of Corpuz et al.~\cite[Section 2.2]{Corpuz25Direc}. For the convergence study against the analytical solution, we construct three meshes, IM0, IM1, and IM2, with successive levels of refinement. The mesh statistics are reported in Table~\ref{tab:meshstatpipe}. A time-step size of $1\times10^{-2}$ is used for all meshes. The cut quadrature of Section~\ref{sec:cut-quad} is used for this study and throughout the paper. Its effectiveness over adaptive subdivision quadrature for internal flows is demonstrated on the pipe flow setup in~\ref{App:AppendixA}. 

\begin{figure}[!t]
    \centering
    \begin{subfigure}{0.495\textwidth}\centering
        \includegraphics[width=\textwidth]{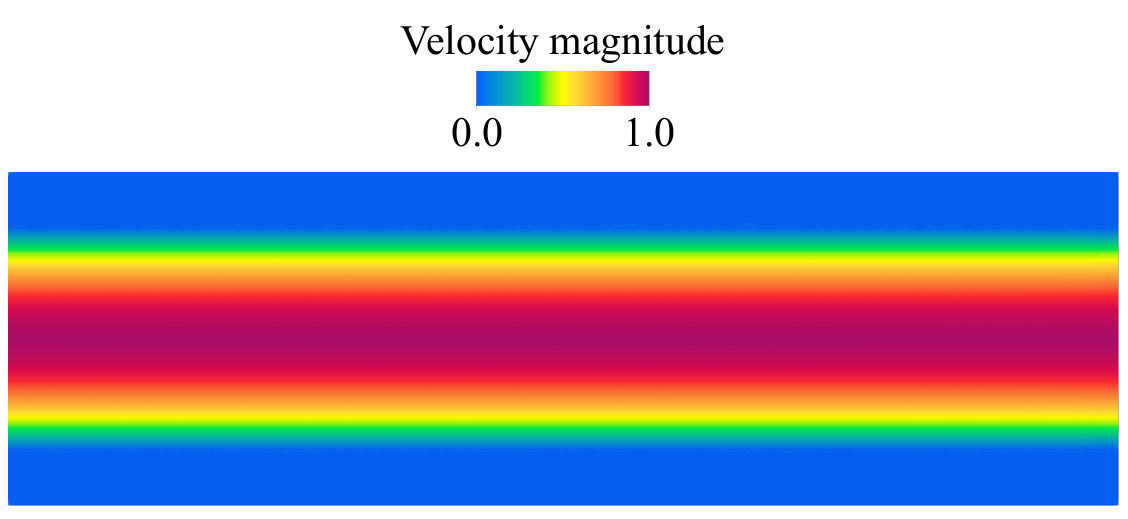}
        \caption{}
        \label{fig:pipe-vel-field}
    \end{subfigure}
    \begin{subfigure}{0.495\textwidth}\centering
        \includegraphics[width=\textwidth]{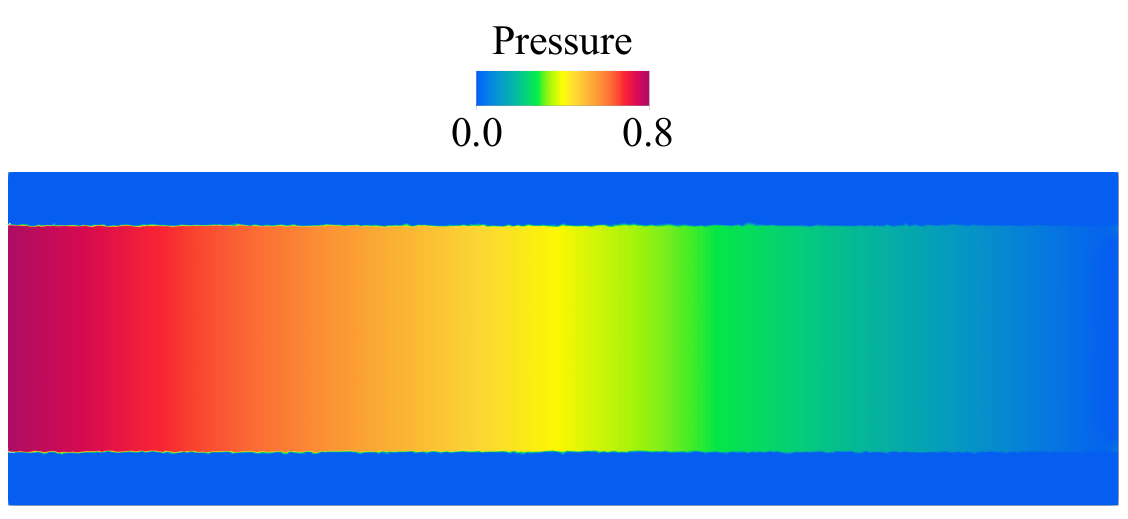}
        \caption{}
        \label{fig:pipe-press-field}
    \end{subfigure}
    \caption{(a) Velocity magnitude and (b) pressure field at the central cross-section for Hagen–Poiseuille pipe flow on IM2. }
    \label{fig:pipe-vel-press-field}
\end{figure}

\begin{figure}[!t]
    \centering
    \begin{subfigure}{0.485\textwidth}\centering
        \includegraphics[width=\textwidth]{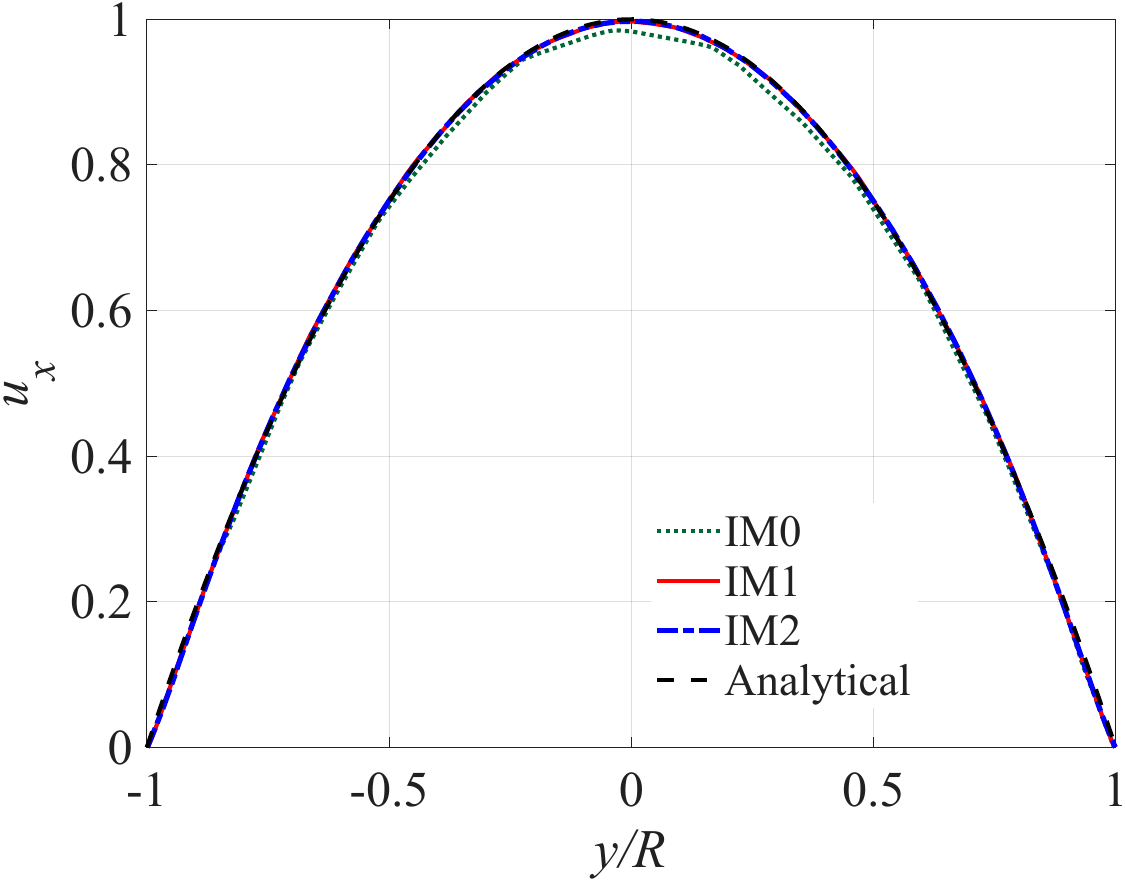}
        \caption{}
        \label{fig:pipe-vel}
    \end{subfigure}
    \hspace{0.01\textwidth}
    \begin{subfigure}{0.485\textwidth}\centering
        \includegraphics[width=\textwidth]{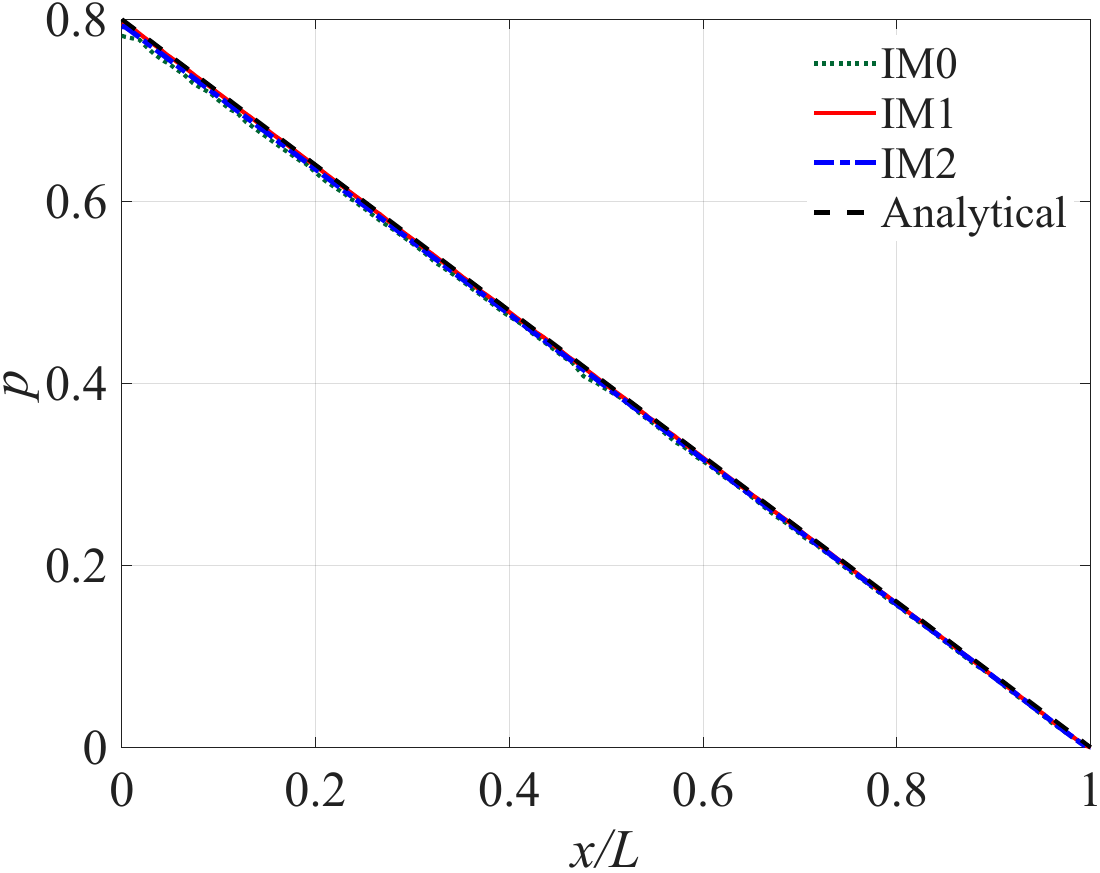}
        \caption{}
        \label{fig:pipe-press}
    \end{subfigure}
    \caption{Convergence of (a) streamwise velocity at the outlet and (b) axial pressure drop compared with the analytical result for Hagen–Poiseuille pipe flow.}
    \label{fig:pipe-vel-press}
\end{figure}

Figure~\ref{fig:pipe-vel-press-field} shows the velocity magnitude and pressure field of the fully developed flow computed on IM2. The fields are visualized on the central cross-section of the entire background domain, comprising both the physical and the fictitious regions. It is evident that the immersed boundary is resolved sharply. The velocity vanishes across the point cloud boundary, the fictitious region carries no spurious flow, and the pressure exhibits the expected linear axial drop within the pipe. The accuracy of the solution is assessed by comparing the outlet velocity profile, the axial pressure drop, and the WSS against the analytical results. Figure~\ref{fig:pipe-vel-press} presents the outlet velocity profile and the pressure drop along the pipe for each mesh level. Both quantities converge strongly under refinement, with only the coarsest mesh, IM0, exhibiting a slight deviation.

\begin{figure}[!t]
    \centering
    \begin{subfigure}{0.485\textwidth}\centering
        \includegraphics[width=\textwidth]{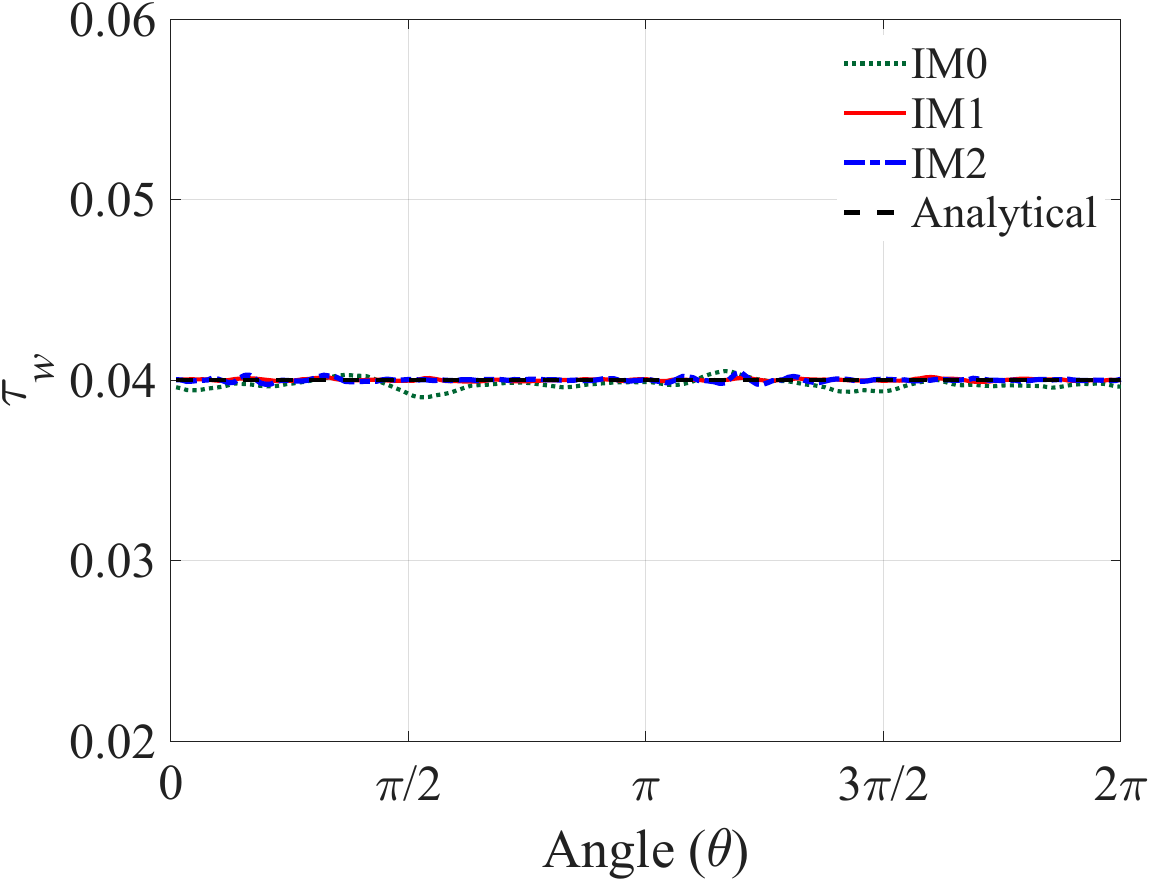}
        \caption{}
        \label{fig:pipe-wss-conv}
    \end{subfigure}
    \hspace{0.01\textwidth}
    \begin{subfigure}{0.485\textwidth}\centering
        \includegraphics[width=\textwidth]{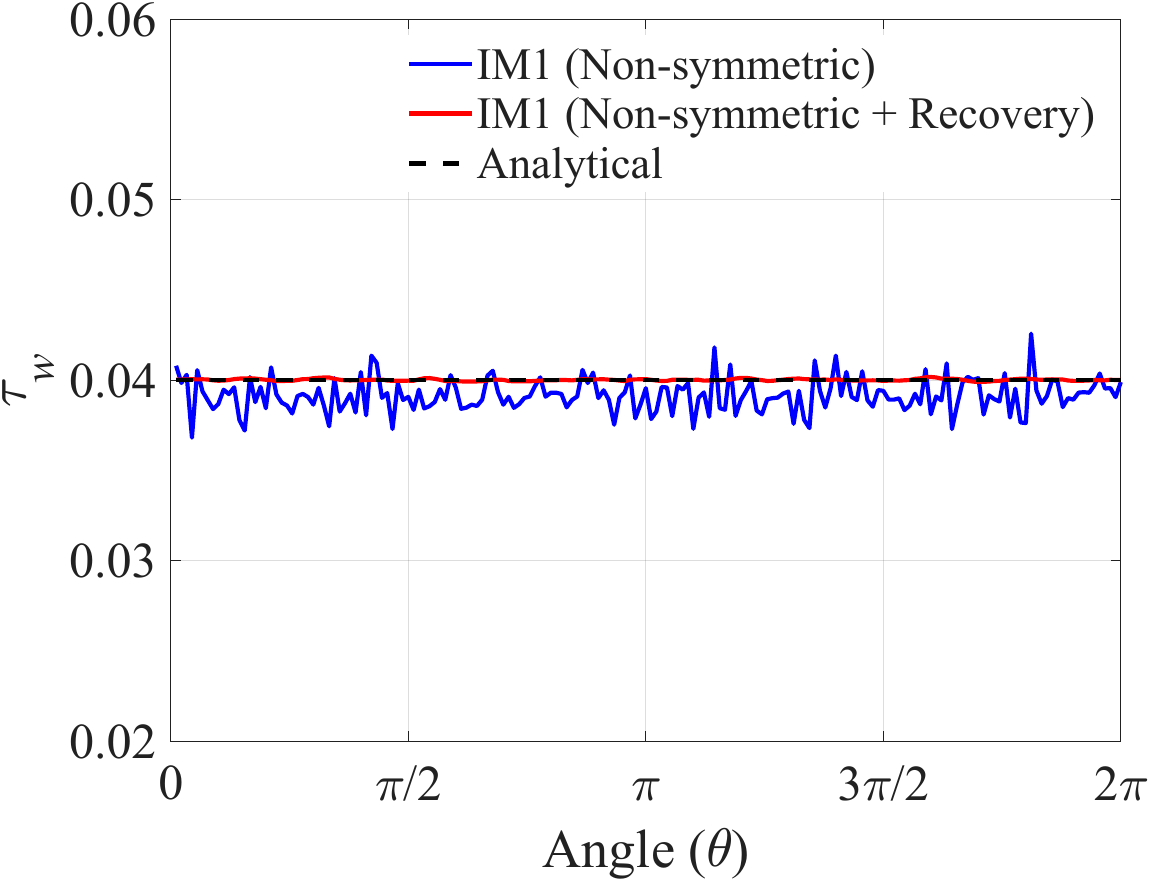}
        \caption{}
        \label{fig:pipe-wss-im1}
    \end{subfigure}
    \caption{(a) Convergence of the recovered WSS around the circumferential cross-section at $x/L = 0.5$ of the pipe compared against the analytical result for Hagen–Poiseuille pipe flow. (b) Comparison of the WSS with and without the proposed WSS recovery method.}
    \label{fig:pipe-wss}
\end{figure}

The recovery procedure of Section~\ref{sec:wss_recovery} is then applied to each mesh, and Figure~\ref{fig:pipe-wss-conv} compares the recovered WSS with the analytical value. The recovered WSS shows excellent convergence, with IM0 again showing only slight variation about the constant analytical value. Figure~\ref{fig:pipe-wss-im1} compares the WSS obtained with and without recovery on the mesh IM1. Without recovery, the WSS oscillates significantly about the exact solution because of the cut elements, but the recovered WSS is both smooth and accurate. We emphasize that the proposed recovery method is not merely a smoothing or averaging procedure. The raw oscillations in Figure~\ref{fig:pipe-wss-im1} are not centered on the exact solution, so any smoothing or averaging of the raw traction would converge to a wrong result. The recovery instead removes the oscillation while yielding a WSS in close agreement with the exact solution. 

\begin{figure}[!t]
    \centering
    \begin{subfigure}{0.485\textwidth}\centering
        \includegraphics[width=\textwidth]{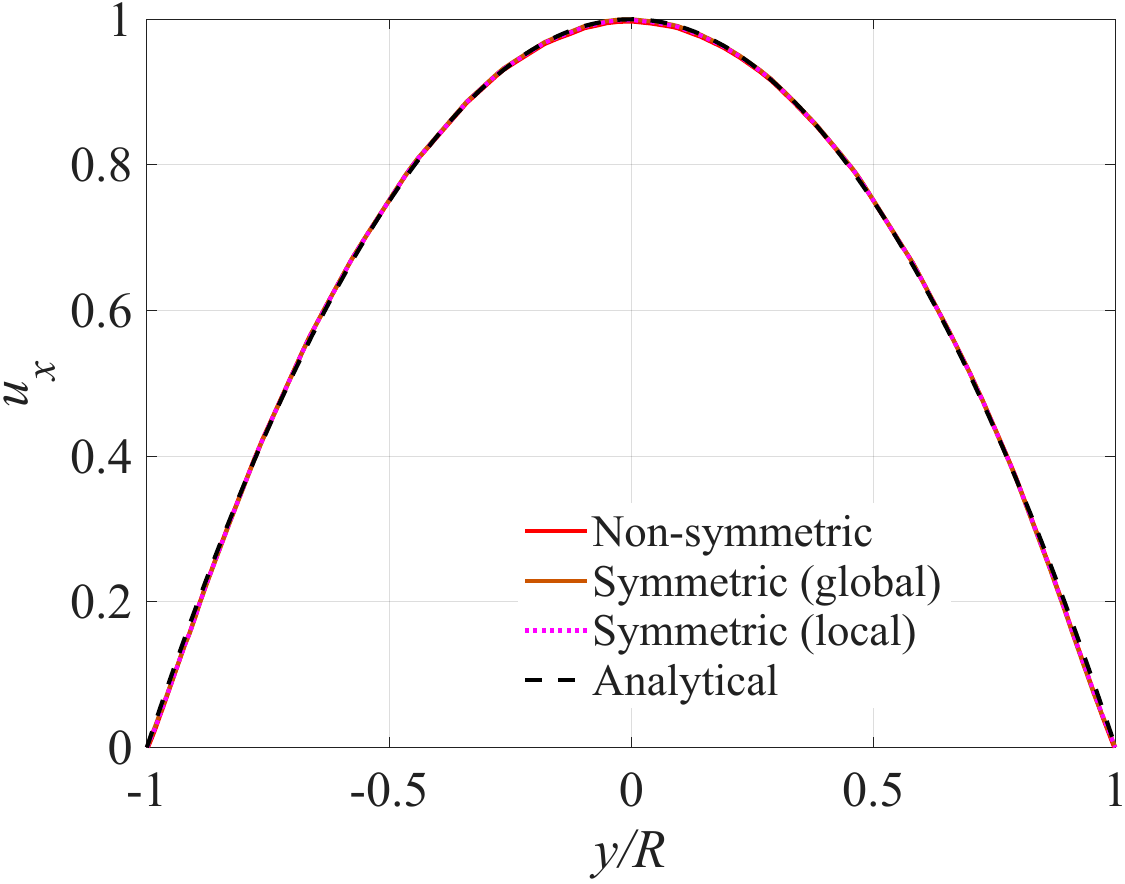}
        \caption{}
        \label{fig:pipe-vel-symm}
    \end{subfigure}
    \hspace{0.01\textwidth}
    \begin{subfigure}{0.485\textwidth}\centering
        \includegraphics[width=\textwidth]{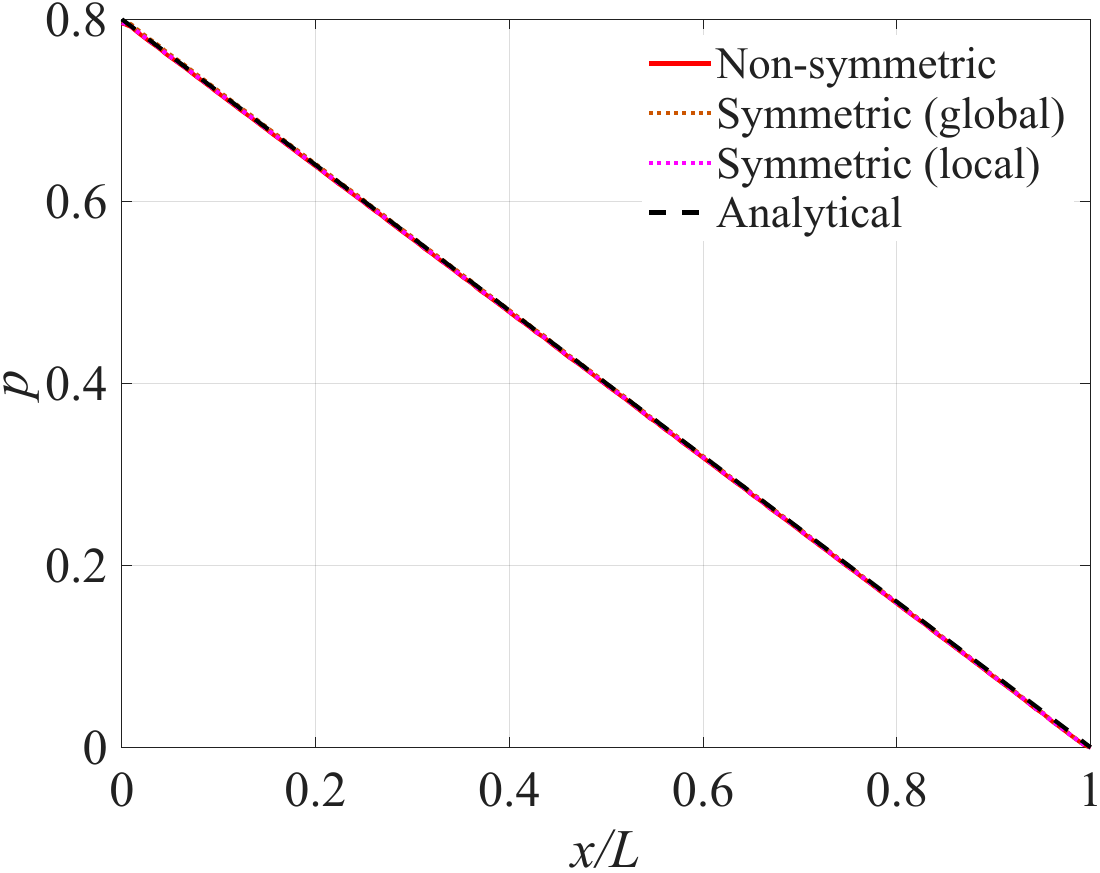}
        \caption{}
        \label{fig:pipe-press-symm}
    \end{subfigure}
    \caption{(a) Streamwise velocity at the outlet and (b) axial pressure drop obtained with the non-symmetric weak BC and with the symmetric weak BC using global and element-local parameter selection, compared with the analytical solution for Hagen–Poiseuille pipe flow on IM1.}
    \label{fig:pipe-symm-vel-press}
\end{figure}
\begin{figure}[!t]
    \centering
    \begin{subfigure}{0.485\textwidth}\centering
        \includegraphics[width=\textwidth]{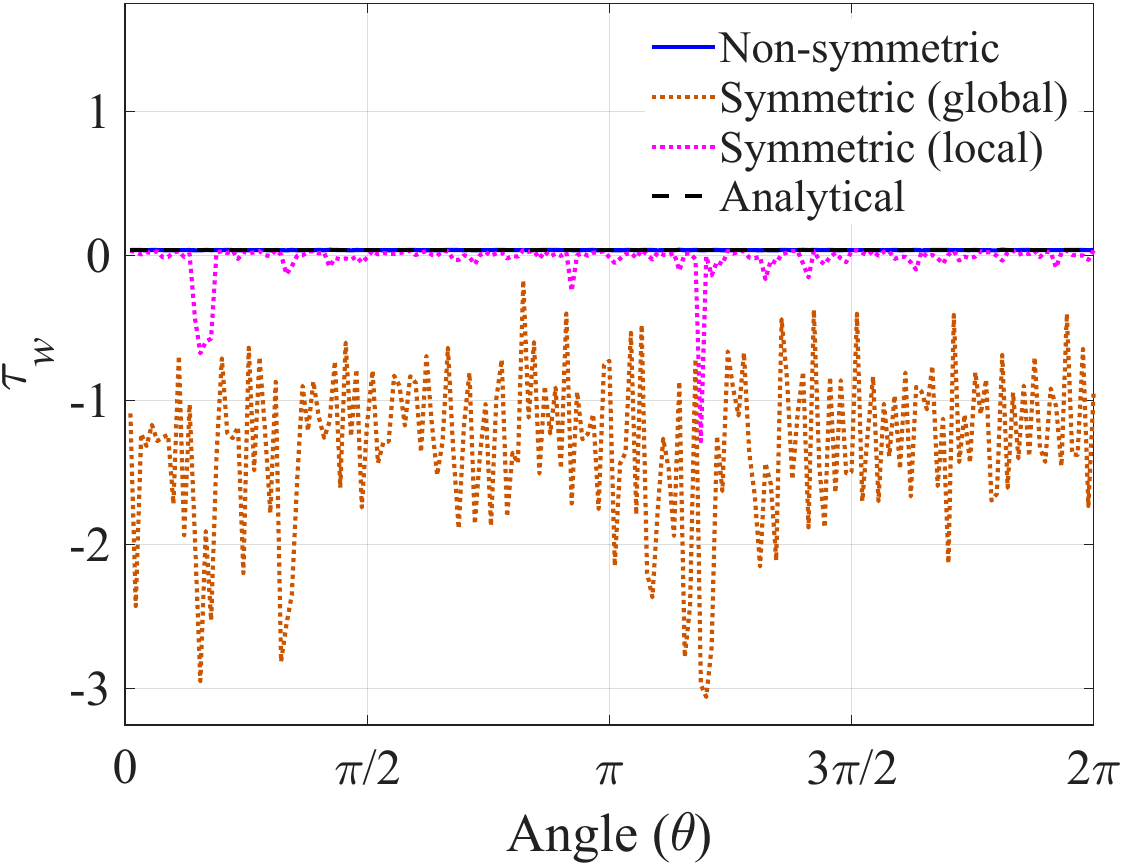}
        \caption{}
    \end{subfigure}
    \hspace{0.01\textwidth}
    \begin{subfigure}{0.4\textwidth}\centering
        \includegraphics[width=\textwidth]{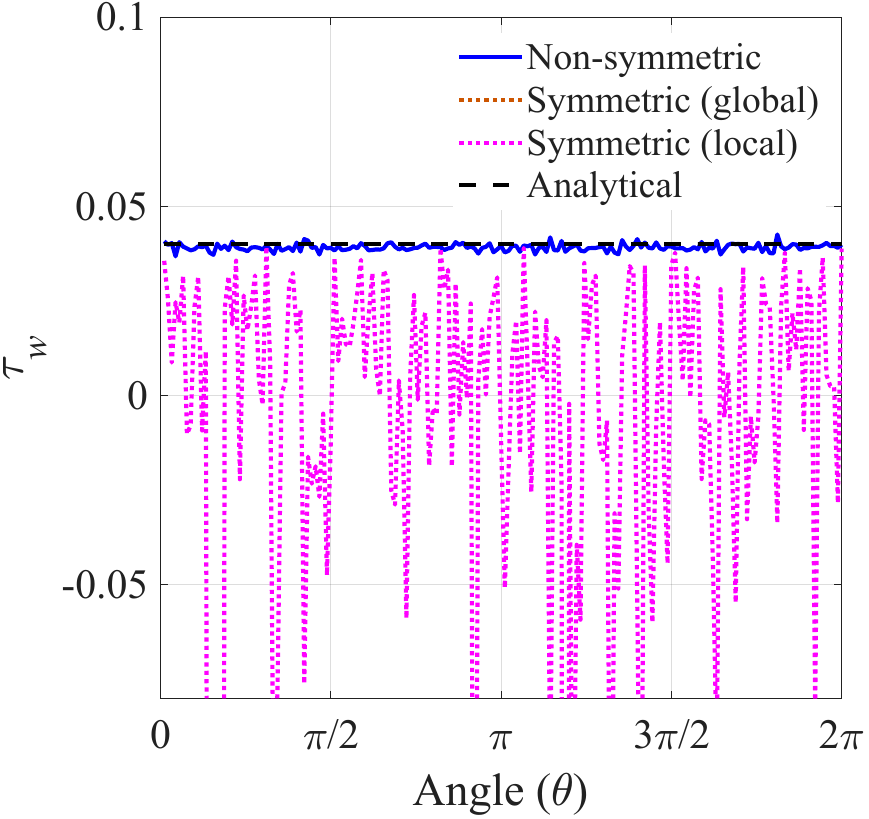}
        \caption{}
    \end{subfigure}
    \caption{Raw WSS (without recovery) around the circumferential cross-section at $x/L = 0.5$, obtained with the non-symmetric weak BC and with the symmetric weak BC using global and element-local parameter selection, compared with the analytical solution for Hagen–Poiseuille pipe flow on IM1. (a) full range, (b) detail near the analytical value.}
    \label{fig:pipe-wss-symm}
\end{figure}

To demonstrate the effectiveness of using non-symmetric Nitsche-based weak BC, we re-compute the flow simulation on IM1 using the symmetric variant ($\gamma = 1$). Two strategies for selecting the stabilization parameters are considered: a sufficiently large global constant that satisfies coercivity everywhere ($\tau^B_\text{NOR} = \tau^B_\text{TAN} = 10^3$)~\cite{Xu16tetra}, and element-local parameters via solution to an eigenvalue problem~\cite{Embar10Impos}. Figure~\ref{fig:pipe-symm-vel-press} shows the outlet velocity profile and the axial pressure drop for both strategies alongside the non-symmetric results, where all three methods produce results close to the analytical solution. However, the raw WSS from each method varies significantly, as shown in Figure~\ref{fig:pipe-wss-symm}. With the global constant parameter, the stabilization term uniformly amplifies the boundary condition residual, and the resulting WSS is dominated by non-physical oscillations and exceeds the analytical result by two orders of magnitude. The element-local estimation performs relatively better in magnitude, but the parameter inherits the sensitivity of the eigenvalue problem to the cut configuration and produces spurious values that affect the local distribution. In contrast, the non-symmetric formulation remains close to the analytical value at this scale even without recovery. These results clearly demonstrate the poor performance of the symmetric Nitsche-based weak BC in predicting accurate WSS in immersed methods, as identified in Section~\ref{sec:weakbc}.

\subsection{Laminar flow around a sphere}
\label{sec:sphere100}
\begin{figure}[!t]
    \centering
    \begin{subfigure}{0.45\textwidth}\centering
        \includegraphics[width=\textwidth]{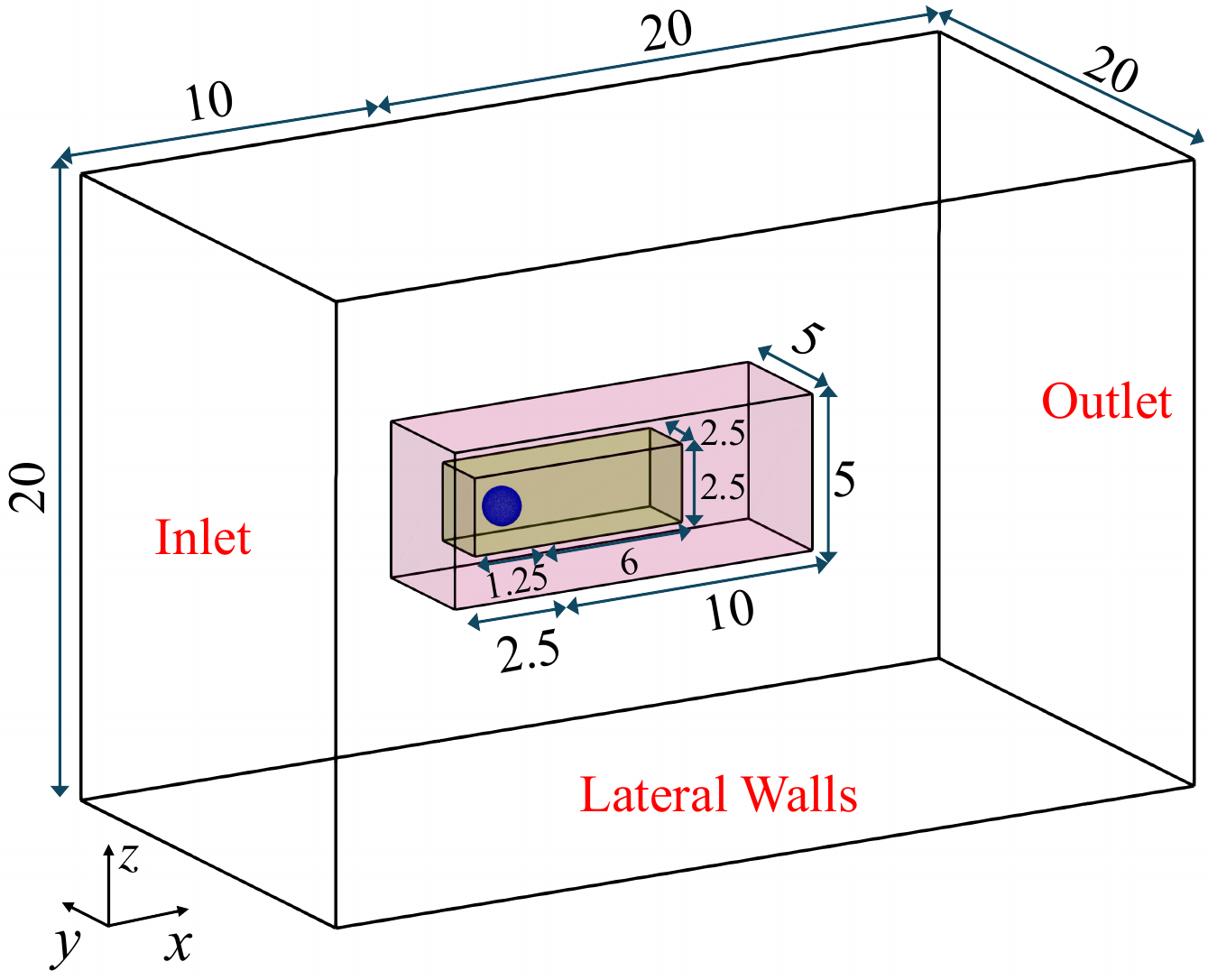}
        \caption{}
        \label{fig:sphere-domain}
    \end{subfigure}
    \hspace{0.01\textwidth}
    \begin{subfigure}{0.525\textwidth}\centering
        \raisebox{0.04\height}{\includegraphics[width=\textwidth]{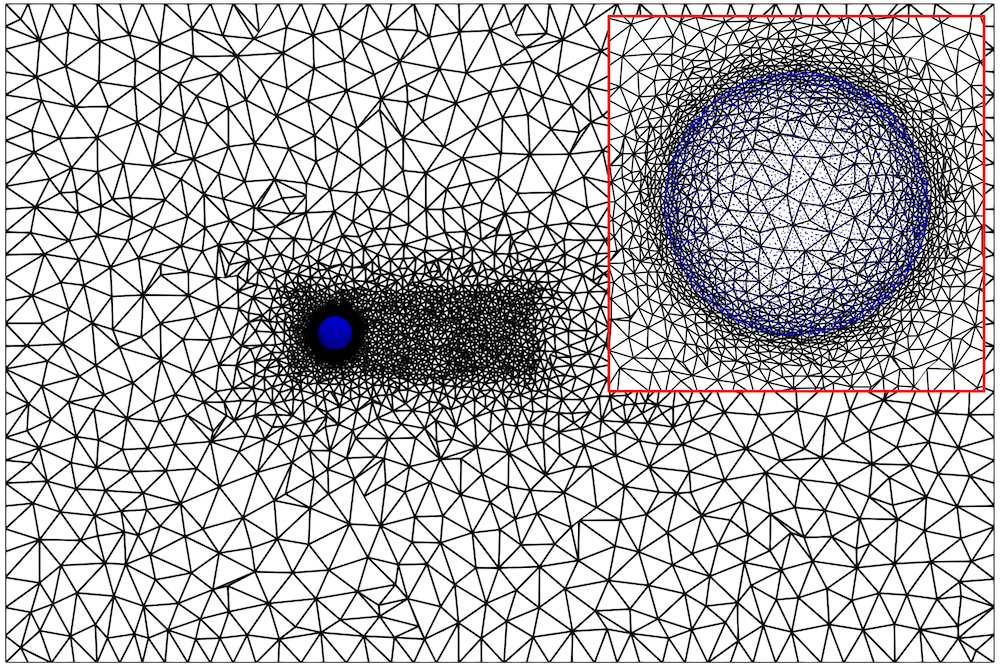}}
        \caption{}
        \label{fig:sphere-im0}
    \end{subfigure}
    \caption{(a) Computational domain and (b) cross-section slice of the background mesh IM0 considered for flow around a sphere, showcasing the domain extent, boundary conditions, and refinement regions.}
    \label{fig:sphere-setup}
\end{figure}

In this study, we consider laminar flow around a sphere at $Re=100$ and compare the solution obtained with the proposed framework against a boundary-fitted reference. The reference solution is computed using the same VMS formulation with symmetric weak BC on a boundary-fitted mesh with prismatic boundary layers around the sphere. For this study, a sphere with diameter $D=1$ is sampled randomly with $10^5$ points over the analytical definition. The resulting point cloud is fully immersed in the computational domain as shown in Figure~\ref{fig:sphere-domain}. The domain comprises two refinement boxes together with near-boundary anisotropic refinement to resolve the boundary layer and the wake behind the sphere geometry. A uniform inflow velocity $U_\infty=1$ is prescribed strongly at the inlet, a traction-free condition at the outlet, and no-penetration conditions on the lateral boundaries. The no-slip boundary condition is imposed on the point cloud using the non-symmetric weak BC (Section~\ref{sec:weakbc}). To produce a laminar flow at $Re=100$, the density and the dynamic viscosity are set to $\rho = 1$ and $\mu = 0.01$, respectively. For the convergence study, the domain is discretized at three successive levels of refinement, producing background meshes IM0, IM1, and IM2. The element sizes and mesh statistics are reported in Table~\ref{tab:meshstatsphere} and the coarsest mesh, IM0, is shown in  Figure~\ref{fig:sphere-im0},  illustrating the near-boundary refinement. The simulations are carried out with a time-step size of $1\times10^{-2}$ until steady state is reached. 

\begin{table}[!t]\centering\small
    \centering
    \captionsetup{justification=centering}
    \newcommand{\tabincell}[2]{\begin{tabular}{@{}#1@{}}#2\end{tabular}}
    \caption{Background mesh statistics for the flow around a sphere at $Re=100$.}
    \begin{tabular}{lccccccr}
        \toprule
        \multirow{2}{*}{Mesh} 
            & \multicolumn{3}{c}{Near wall} 
            & \multirow{2}{*}{\tabincell{c}{Inner \\ refinement box}}
            & \multirow{2}{*}{\tabincell{c}{Outer \\ refinement box}}
            & \multirow{2}{*}{Domain} 
            & \multirow{2}{*}{\tabincell{c}{Number of \\ active elements}} \\
        \cmidrule(lr){2-4}
            & $h_\mathrm{NOR}$
            & $h_\mathrm{TAN}$
            & \tabincell{c}{BL height} 
            & & \\
        \midrule
        IM0 & $0.02$  & $0.04$ & $0.04$ & $0.2$  & $0.8/\sqrt{2}$ & $1.2$ & $378{,}755$  \\
        IM1 & $0.01$  & $0.02$ & $0.02$ & $0.1$  & $0.4/\sqrt{2}$ & $1.0$ & $1{,}461{,}858$ \\
        IM2 & $0.005$ & $0.01$ & $0.01$ & $0.05$ & $0.2/\sqrt{2}$ & $0.8$ & $7{,}032{,}884$ \\
        \bottomrule
    \end{tabular}
    \label{tab:meshstatsphere}
\end{table}

\begin{figure}[!t]
    \centering
    \begin{subfigure}{0.575\textwidth}\centering
        \includegraphics[width=\textwidth]{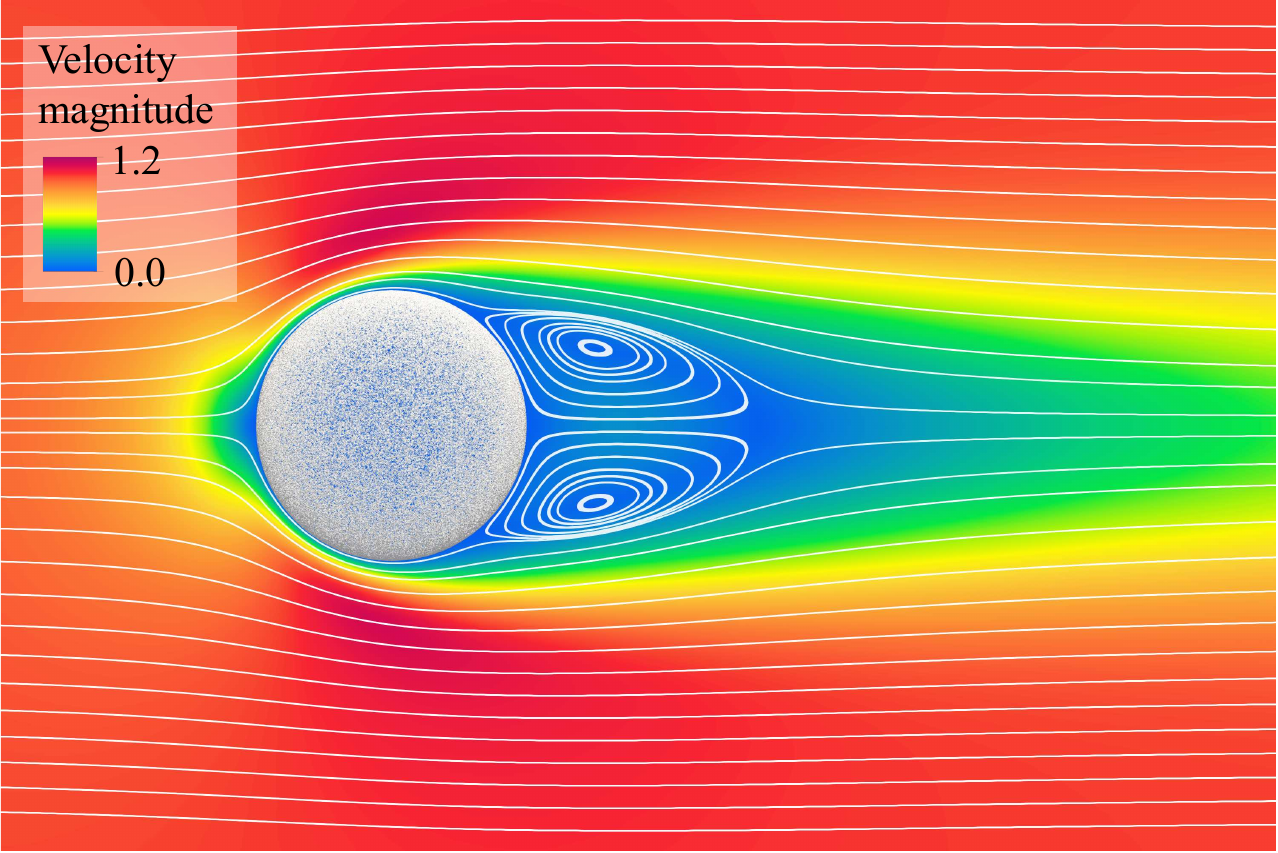}
        \caption{}
        \label{fig:sphere-re100-vel}
    \end{subfigure}
    \hspace{0.01\textwidth}
    \begin{subfigure}{0.4\textwidth}\centering
        \includegraphics[width=\textwidth]{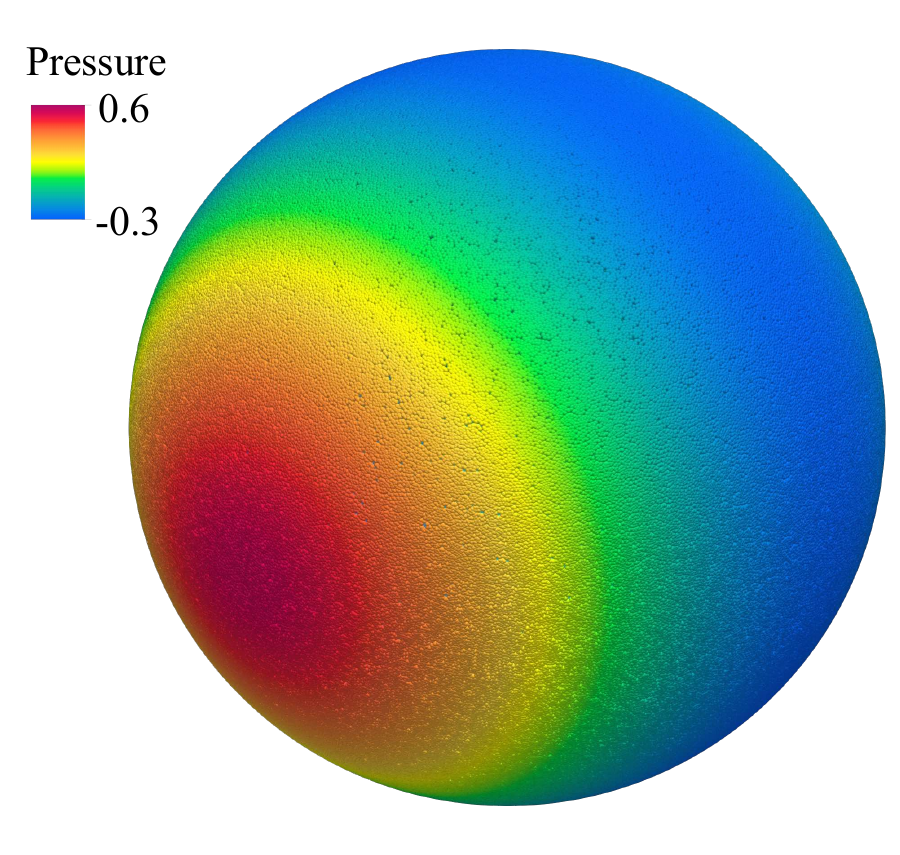}
        \caption{}
        \label{fig:sphere-re100-press}
    \end{subfigure}
    \caption{(a) Velocity magnitude on the central cross-section and (b) pressure field evaluated on the point cloud for flow around a sphere at $Re=100$.}
    \label{fig:sphere-re100-vel-press}
\end{figure}

\begin{figure}[!t]
    \centering
    \begin{subfigure}{0.485\textwidth}\centering
        \includegraphics[width=\textwidth]{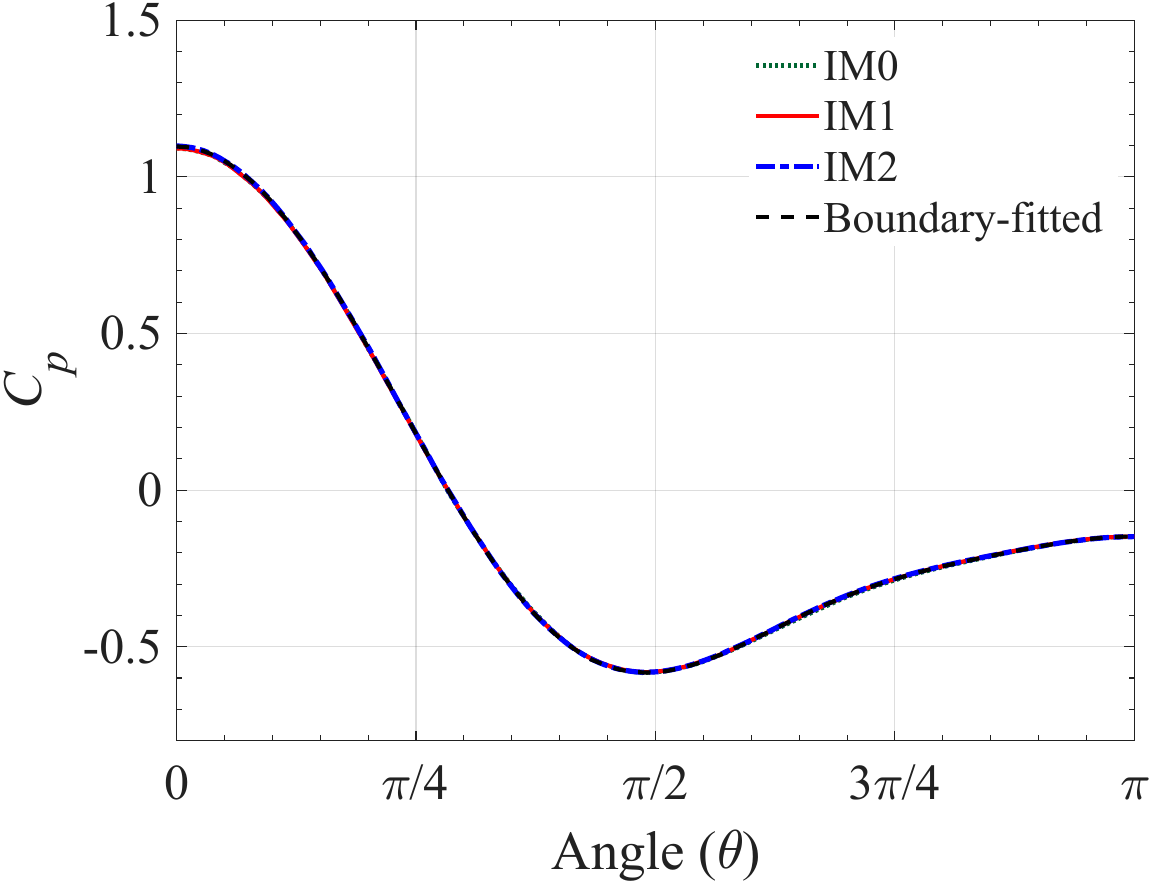}
        \caption{}
        \label{fig:sphere-re100-cp}
    \end{subfigure}
    \hspace{0.01\textwidth}
    \begin{subfigure}{0.485\textwidth}\centering
        \includegraphics[width=\textwidth]{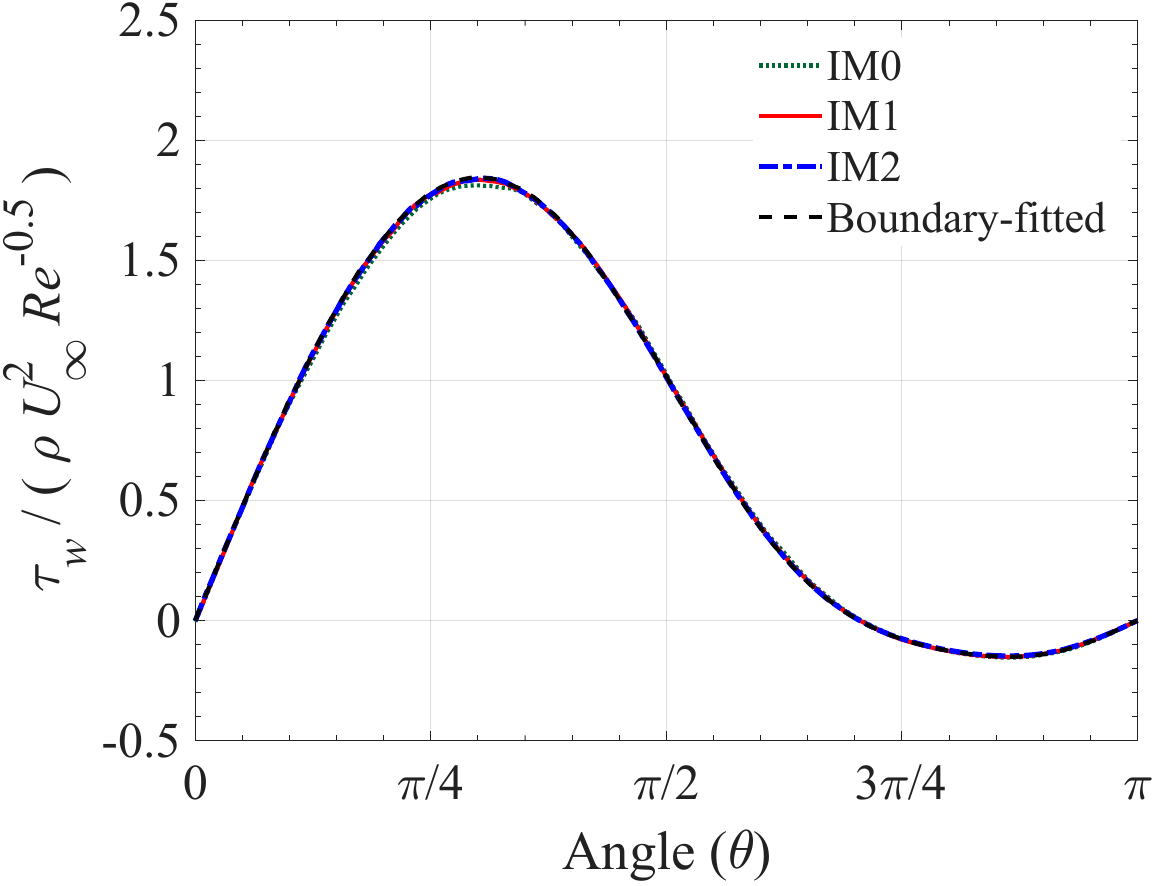}
        \caption{}
        \label{fig:sphere-re100-cf}
    \end{subfigure}
    \caption{Convergence of (a) pressure coefficient and (b) skin friction evaluated on the upper crown line of the sphere for flow around a sphere at $Re=100$.}
    \label{fig:sphere-re100-cp-cf}
\end{figure}

Figure~\ref{fig:sphere-re100-vel} shows the velocity magnitude contours and streamlines on the central cross-section of the domain. The field near the point cloud boundary is resolved sharply, and the axisymmetric recirculation bubble behind the sphere is captured. The pressure field evaluated on the point cloud is presented in Figure~\ref{fig:sphere-re100-press}, showing a smooth distribution without any numerical artifacts or spurious oscillations. To demonstrate the accuracy of the solution, we compute the drag coefficient ($C_D$) by integrating the variationally consistent traction over the point cloud and the pressure coefficient ($C_p$) along the upper crown line of the sphere. Table~\ref{tab:spherere100_drag} reports the $C_D$ and Figure~\ref{fig:sphere-re100-cp} shows the $C_p$ obtained for each background mesh. Both quantities converge with mesh refinement and agree closely with the reference results. Next, the proposed recovery method is applied to extract the WSS on the point cloud, from which the skin friction $\left(\tau_w / ( \rho U_\infty^2 Re^{-0.5}) \right)$ is evaluated along the upper crown line of the sphere for each mesh, as shown in Figure~\ref{fig:sphere-re100-cf}. The recovered skin friction also exhibits excellent mesh convergence and agreement with the boundary-fitted reference. 

\begin{table}[!t]\centering\small
    \centering
    \captionsetup{justification=centering}
    \newcommand{\tabincell}[2]{\begin{tabular}{@{}#1@{}}#2\end{tabular}}
    \caption{Drag coefficient obtained for the flow around a sphere at $Re= 100$.}
    \begin{tabular}{ll}
        \toprule
         \tabincell{l}{Mesh} & \tabincell{l}{$C_D$} \\
         \midrule
         \tabincell{l}{IM0} & \tabincell{l}{1.102} \\
         \tabincell{l}{IM1} & \tabincell{l}{1.095} \\
         \tabincell{l}{IM2} & \tabincell{l}{1.094} \\
        \midrule
        \tabincell{l}{Boundary-fitted} & \tabincell{l}{1.094}  \\
        \bottomrule
    \end{tabular}
    \label{tab:spherere100_drag}
\end{table}

\begin{figure}[!t]
    \centering
    \includegraphics[width=\textwidth]{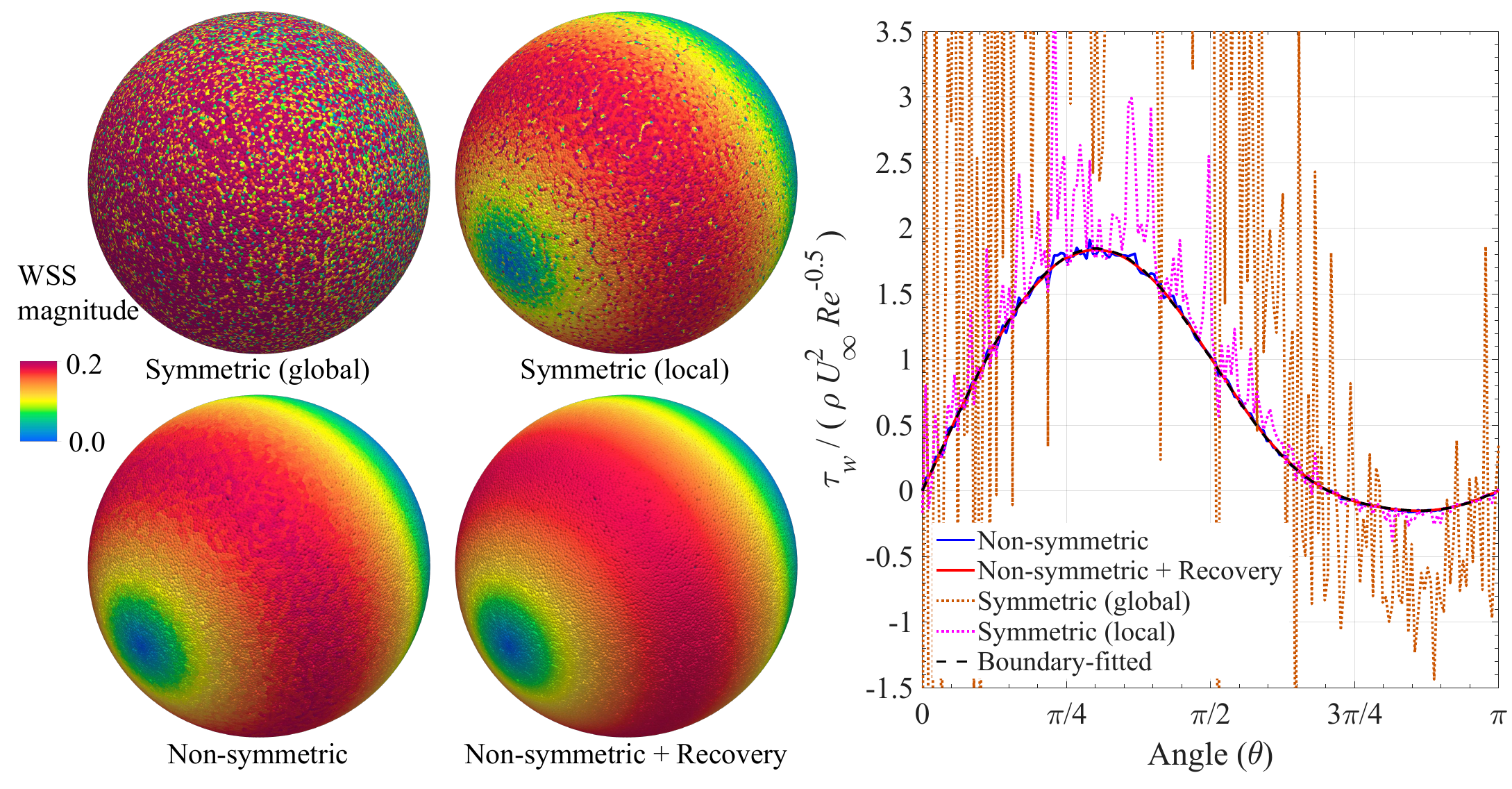}
    \caption{Comparison of WSS magnitude evaluated on the point cloud (left) and skin friction plotted on the upper crown line (right) of the sphere obtained with different weak BC formulations for flow around a sphere at $Re=100$.}
    \label{fig:sphere-re100-wss-compare}
\end{figure}

Similar to the pipe flow, we recompute the flow around the sphere at $Re=100$ on IM1 using the symmetric weak BC, and compare against the non-symmetric variant with and without recovery. Both global constant and element-local, eigenvalue-based parameter selection are again considered for the symmetric formulation. Figure~\ref{fig:sphere-re100-wss-compare} shows the resulting WSS magnitude and skin-friction coefficient for each formulation. Notably, the symmetric weak BC is significantly inferior to the non-symmetric formulation in both cases: the global parameter selection produces WSS dominated by non-physical noise, while the element-local estimation suffers from spurious oscillation arising from cut-element sensitivity. Meanwhile, the non-symmetric formulation yields a substantially smoother and more accurate WSS estimation even without recovery, and the proposed recovery method further improves it, producing a WSS in close agreement with the boundary-fitted reference. These results confirm the effectiveness of the proposed method in extracting accurate WSS on a fully three-dimensional problem with a curved boundary. 

\section{Applications}
\label{sec:applications}
This section demonstrates the proposed framework on higher Reynolds number flows and on complex, scan-derived geometry. We consider two cases: turbulent flow past a sphere at \mbox{$Re=3700$}, and flow in a patient-specific aorta, where accurate WSS is of direct clinical relevance. Both cases are compared against reference results. 

\subsection{Turbulent flow past a sphere at $Re=3700$}
\label{sec:sphere3700}
Flow past a sphere at $Re=3700$ lies in the subcritical regime, where the boundary layer undergoes laminar separation, producing a detached shear layer that transitions to a turbulent wake with irregular vortex shedding. This behavior makes the case a demanding benchmark for numerical methods, and it has accordingly been studied extensively~\cite{Rodriguez11Direc, Yun06Vorti, Bazilevs14Compu}. The case is also a stringent test for WSS estimation. The laminar boundary layer is thin at this Reynolds number, forming a very steep wall-normal velocity gradient confined to the region where the cut elements do not align. The WSS rises sharply over the upstream face and remains small and fluctuating in the separated region. Accurate WSS estimation therefore requires both faithful near-wall treatment and a boundary stress evaluation free of discretization artifacts, which the proposed framework provides through the wall-model-based parameter selection and the WSS recovery, respectively. 

We utilize the same computational setup as in Section~\ref{sec:sphere100} where a randomly sampled point cloud of a sphere is immersed into a domain, with inflow velocity $U_\infty = 1$, density $\rho=1$, and dynamic viscosity $\mu = 1/Re$. Similarly, the domain is discretized into tetrahedral elements with refinement boxes and near-boundary anisotropic refinements. The element size and mesh statistics considered for this study are reported in Table~\ref{tab:meshstatsphere3700}. A time-step size of $1\times 10^{-3}$ is used in this study.

\begin{table}[!t]\centering\small
    \centering
    \captionsetup{justification=centering}
    \newcommand{\tabincell}[2]{\begin{tabular}{@{}#1@{}}#2\end{tabular}}
    \caption{Background mesh statistics for the flow past a sphere at $Re=3700$.}
    \begin{tabular}{ccccccr}
        \toprule
        \multicolumn{3}{c}{Near wall} 
        & \multirow{2}{*}{\tabincell{c}{Inner \\ refinement box}}
        & \multirow{2}{*}{\tabincell{c}{Outer \\ refinement box}}
        & \multirow{2}{*}{Domain} 
        & \multirow{2}{*}{\tabincell{c}{Number of \\ active elements}} \\
        \cmidrule(lr){1-3}
            $h_\mathrm{NOR}$
            & $h_\mathrm{TAN}$
            & \tabincell{c}{BL height} 
            & & \\
        \midrule
        $0.002$  & $0.01$ & $0.01$  & $0.04$ & $0.16/\sqrt{2}$ & $0.8$ & $13{,}350{,}382$ \\
        \bottomrule
    \end{tabular}
    \label{tab:meshstatsphere3700}
\end{table}

\begin{figure}[!t]
    \centering
    \includegraphics[width=\textwidth]{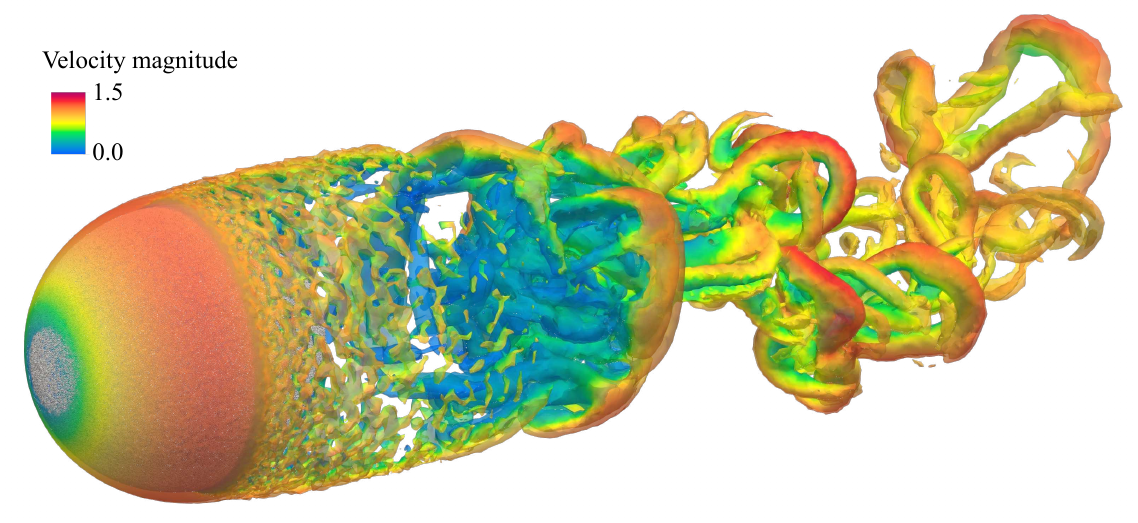}
    \caption{$Q$-criterion of flow past a sphere at $Re=3700$ colored with velocity magnitude. }
    \label{fig:sphere-re3700-qcrit}
\end{figure}

\begin{figure}[!t]
    \centering
    \begin{subfigure}{0.485\textwidth}\centering
        \includegraphics[width=\textwidth]{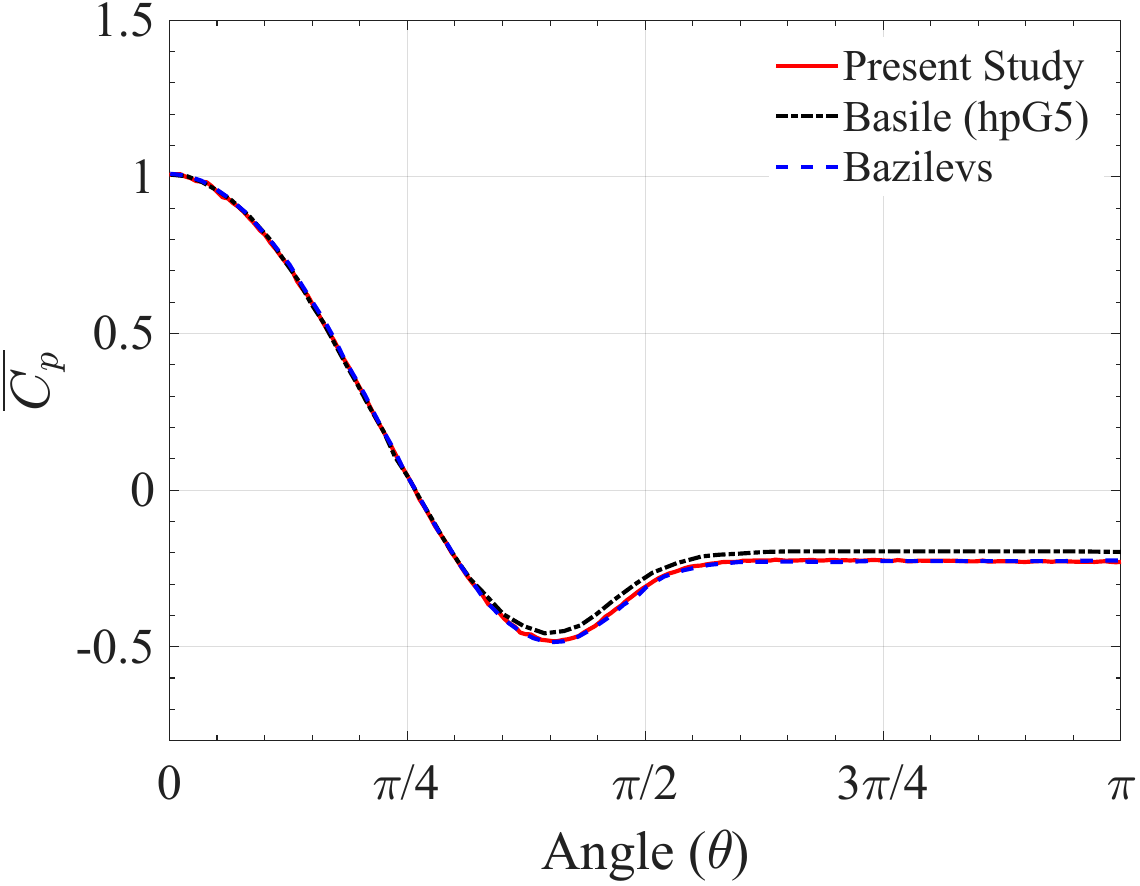}
        \caption{}
        \label{fig:sphere-re3700-press}
    \end{subfigure}
    \hspace{0.01\textwidth}
    \begin{subfigure}{0.485\textwidth}\centering
        \includegraphics[width=\textwidth]{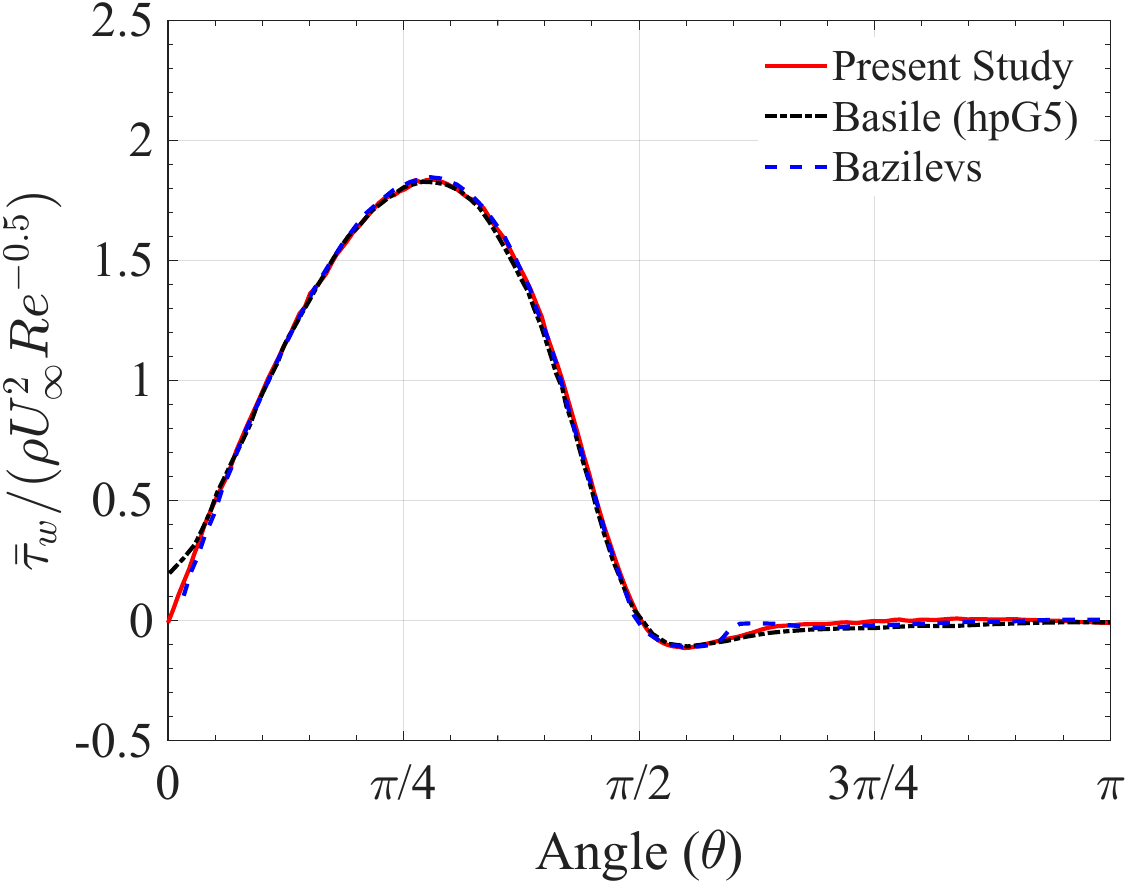}
        \caption{}
        \label{fig:sphere-re3700-wss}
    \end{subfigure}
    \caption{ Time-averaged (a) pressure coefficient and (b) skin-friction for flow past a sphere at $Re=3700$.}
    \label{fig:sphere-re3700-press-wss}
\end{figure}

Figure~\ref{fig:sphere-re3700-qcrit} visualizes the instantaneous wake by the isosurfaces of the $Q$-criterion colored by the velocity magnitude. The framework captures the characteristic features of this flow regime, such as a laminar boundary layer upstream of the sphere, followed by separation leading to a turbulent wake behind the sphere. Figure~\ref{fig:sphere-re3700-press-wss} presents quantitative comparisons of the time-averaged pressure and skin friction coefficients along the upper crown line of the sphere, compared with reference data from the literature~\cite{Bazilevs14Compu, Basile24Scale}. The pressure coefficient (Figure~\ref{fig:sphere-re3700-press}) is in close agreement with the reference throughout the surface, including the front stagnation point, the pressure gradient region, and the base pressure plateau behind the separation. The skin friction coefficient (Figure~\ref{fig:sphere-re3700-wss}), obtained with the proposed recovery method, likewise follows the reference closely. The framework predicts the peak skin friction at $\theta_p = 50.3^\circ$, within the range of $48^\circ$ to $51^\circ$ reported in literature~\cite{Rodriguez11Direc, Bazilevs14Compu}, and a separation angle $\varphi_s = 90.5^\circ$, in good agreement with the value of $89.4^\circ$ obtained both by the DNS~\cite{Rodriguez11Direc} and LES~\cite{Bazilevs14Compu}, and falls within the range reported in Ref.~\cite{Basile24Scale}.

\begin{figure}[!t]
    \centering
    \begin{subfigure}{0.425\textwidth}\centering
        \includegraphics[width=\textwidth]{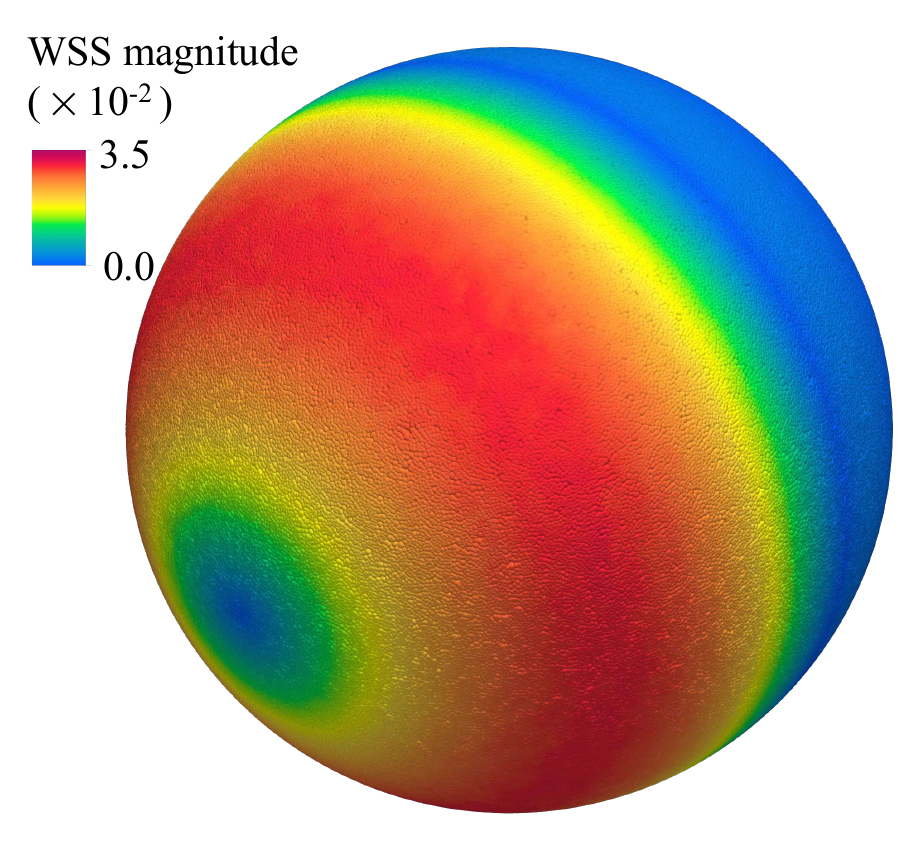}
        \caption{}
        \label{fig:sphere-re3700-wss-pc}
    \end{subfigure}
    \hspace{0.01\textwidth}
    \begin{subfigure}{0.485\textwidth}\centering
        \includegraphics[width=\textwidth]{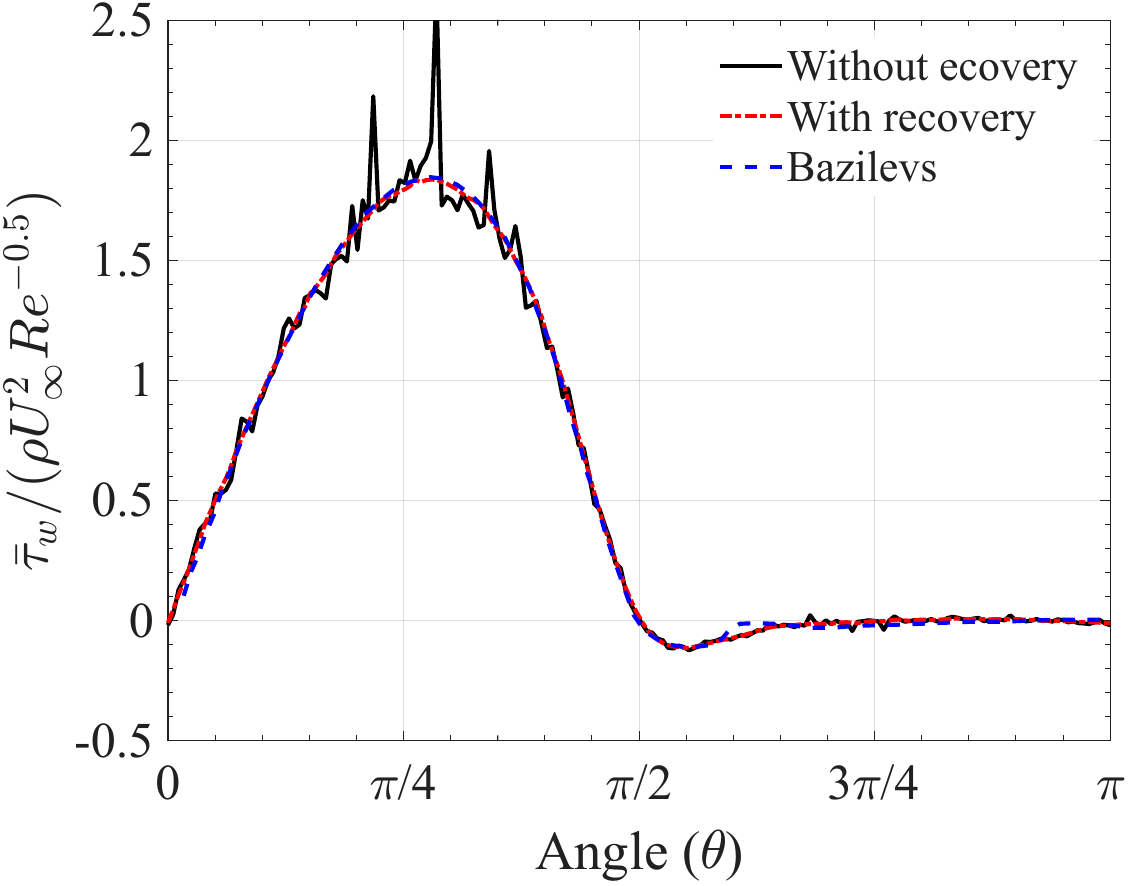}
        \caption{}
        \label{fig:sphere-re3700-pre-post-wss}
    \end{subfigure}
    \caption{ Time-averaged (a) WSS on the point cloud and (b) comparison of WSS with and without recovery for flow past a sphere at $Re=3700$.}
    \label{fig:sphere-re3700-wss-comp}
\end{figure}

The time-averaged recovered WSS distribution on the point cloud is shown in Figure~\ref{fig:sphere-re3700-wss-pc}. The field is smooth over the entire surface and exhibits the expected axisymmetric distribution. Figure~\ref{fig:sphere-re3700-pre-post-wss} compares the time-averaged WSS before and after recovery along the upper crown line of the sphere. Without the recovery, the WSS exhibits pronounced oscillations, particularly near the peak region, reflecting the cut-element irregularity. These oscillations persist even after time averaging, confirming that they are a spatial artifact of the discretization rather than a temporal fluctuation. The recovered WSS removes them while preserving the underlying distribution, yielding a field from which the separation angle and the skin friction distribution can be extracted accurately. 

\subsection{Patient-specific aorta}
\label{sec:aorta}
We now consider the hemodynamics in a patient-specific aorta from a healthy 11-year-old male, previously studied in Refs.~\cite{LaDisa11Compu, Africa24lifex}. The model is obtained from the Vascular Model Repository (VMR)~\cite{Wilson13Vascu}, where an MRI scan was manually segmented using SimVascular~\cite{Updegrove17SimVa} to generate a boundary-fitted mesh for flow analysis. Deep learning-based auto-segmentation methods allow automatic extraction of a point cloud from medical images~\cite{An25Hiera, Du25AIpow}, which can then be analyzed directly with point cloud-based CFD, as demonstrated by Corpuz et al.~\cite{Corpuz25Direc}. Here, however, we deliberately adopt the VMR model, representing its surface by a point cloud and analyzing it with the proposed framework, so that the relevant quantities of interest can be compared one-to-one against boundary-fitted results computed on the same geometry.

\subsubsection{VMR model and computational setup}
\label{sec:aorta-setup}
\begin{figure}[!t]
    \centering
    \includegraphics[width=\textwidth]{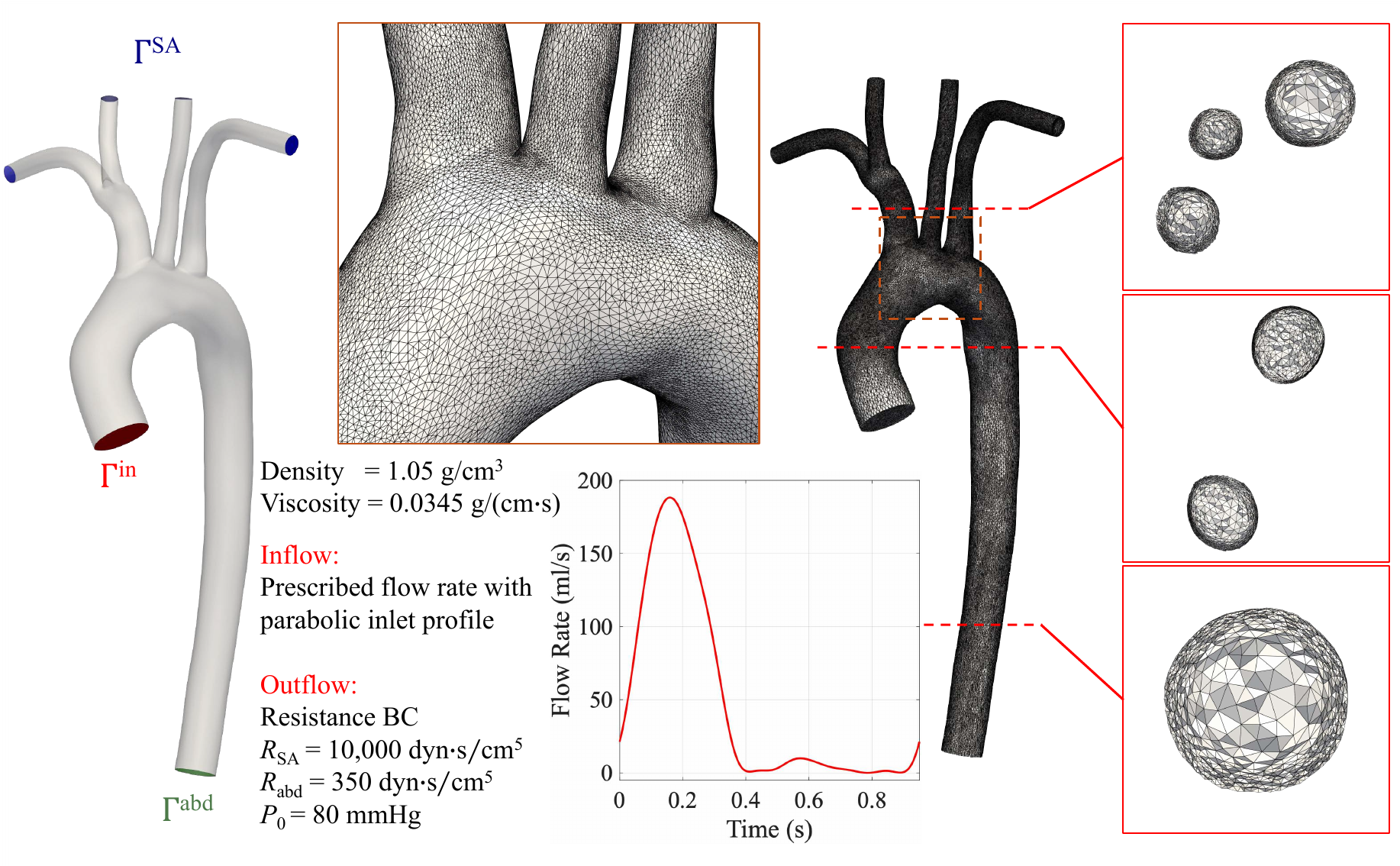}
    \caption{Patient-specific aorta, computational setup for full cardiac cycle, and boundary-fitted mesh with adaptive refinements obtained from the Vascular Model Repository~\cite{Wilson13Vascu}. Cross-section of the mesh at different heights shows boundary layer refinements.}
    \label{fig:aorta-vmr-setup}
\end{figure}

\begin{figure}[!t]
    \centering
    \includegraphics[width=\textwidth]{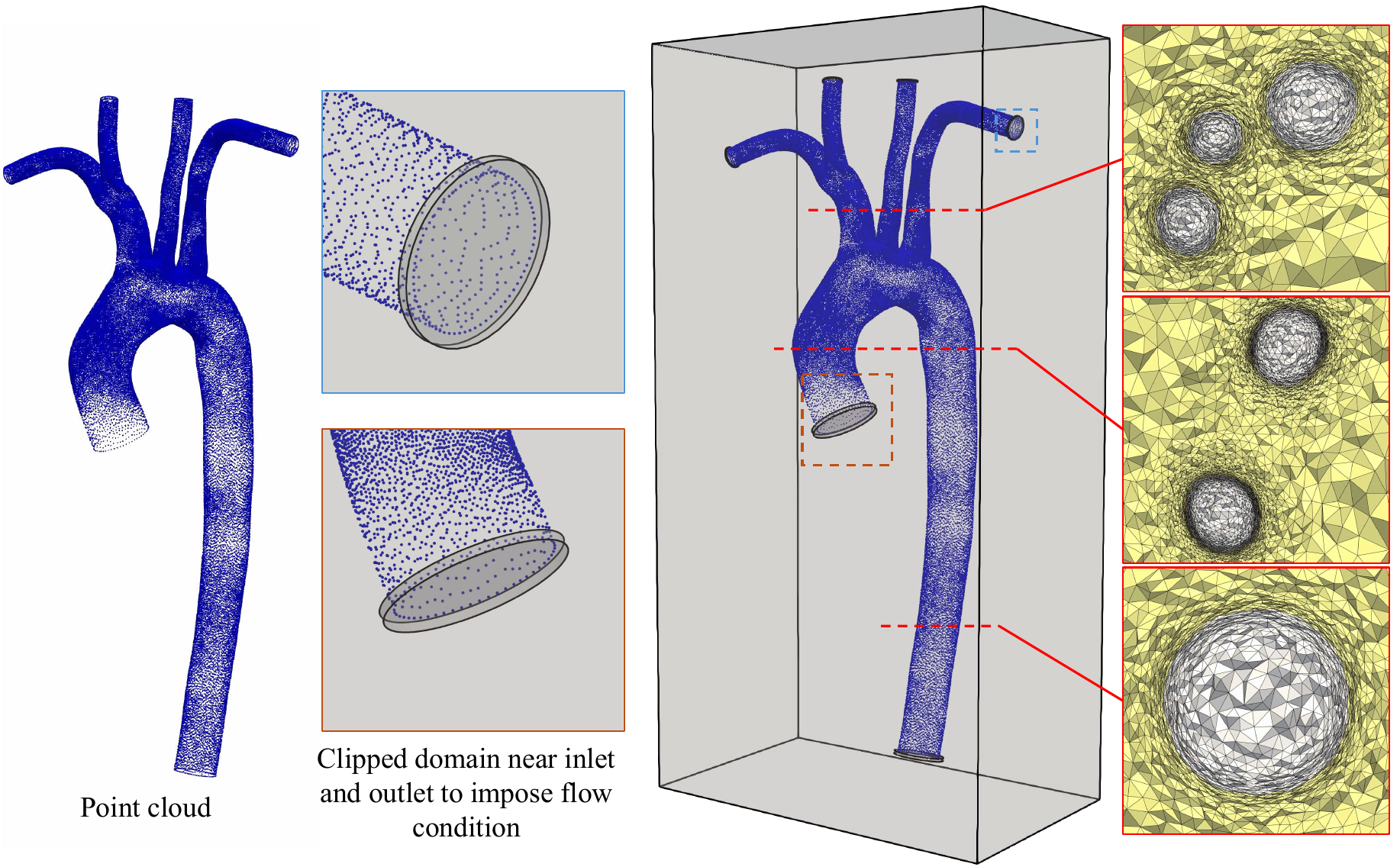}
    \caption{Point cloud generated from sampling the surface mesh vertices of the VMR model and background domain constructed with clipped regions near the inlet and outlet to impose the flow conditions exactly. Cross-section of the background mesh at different heights shows similar discretization to the boundary-fitted mesh. Active elements are shaded white, and inactive elements are shaded yellow.}
    \label{fig:aorta-im-setup}
\end{figure}

Figure~\ref{fig:aorta-vmr-setup} presents the VMR model and computational setup used to simulate a full cardiac cycle. The model features one inlet $\Gamma^\text{in}$ at the aortic root and five outlets, corresponding to the four supra-aortic arteries, collectively denoted $\Gamma^\text{SA}$, and one abdominal aorta $\Gamma^\text{abd}$. The boundary-fitted mesh of the VMR model consists of $2{,}915{,}690$ tetrahedral elements. The mesh was refined near the wall to resolve the boundary layer, with additional surface refinement near the aortic arch to capture the flow along the curvature and branching of the supra-aortic vessels. This refinement was obtained with an adaptive strategy reported in Ref.~\cite{LaDisa11Compu}. A Dirichlet boundary condition is imposed strongly at the inlet $\Gamma^\text{in}$, with a parabolic velocity profile and pulsatile flow rate over time. The parabolic profile is defined on the largest circle inscribed in $\Gamma^\text{in}$. The flow waveform is obtained from the VMR, which was derived from experimental measurements reported in Ref.~\cite{LaDisa11Compu}. It is characterized by a strong ejection phase during ventricular systole, followed by a diastolic phase with nearly zero flow rate. To impose physiological flow conditions at the outlets, we prescribe resistance boundary conditions by adding the following terms to the variational form \eqref{eq:weak_form}:
\begin{align}
\nonumber \sum_a \left( 
- \beta \int_{\Gamma^a} \mathbf{w}^h \cdot \rho \, \left\{ \mathbf{u}^h \cdot \mathbf{n} \right\}_{-} \, \mathbf{u}^h \, d\Gamma \, 
+ \int_{\Gamma^a} \mathbf{w}^h \cdot \mathbf{n} \left( p_0 + R_a \int_{\Gamma^a} \mathbf{u}^h \cdot \mathbf{n} \, d\Gamma \right) \, d\Gamma
\right) \text{ ,}
\end{align}
where $\Gamma^a$ denotes the outflow portions of the domain, $\beta \in [0,1]$ is a dimensionless scalar, $\{\cdot\}_{-}$ takes the negative part of its argument, and $p_0$ and $R_a$ are the base pressure and the resistance prescribed on $\Gamma^a$, respectively. The first term is a backflow stabilization that controls the unstable convective influx at the outflow boundary and prevents solver divergence~\cite{Bazilevs09Patie}. In this work, we set $\beta = 0.5$ following the approach of Esmaily-Moghadam et al.~\cite{Esmaily11compa}. The second term enforces the resistance boundary condition at the outlets. We set $R_\text{SA} = 10^4 \text{ dyn} \cdot \text{s} / \text{cm}^5 $ on the supra-aortic outlets $\Gamma^\text{SA}$ and $R_\text{abd} = 3.5 \times 10^2 \text{ dyn} \cdot \text{s} / \text{cm}^5 $ on the abdominal outlet $\Gamma^\text{abd}$, with $p_0= 80 \text{ mmHg}$. The density and dynamic viscosity are set to $\rho = 1.05 \text{ g}/\text{cm}^3$ and $\mu = 3.45 \times 10^{-2} \text{ g}/ (\text{cm} \cdot \text{s})$, respectively. The simulation is carried out with a time-step size of $1\times 10^{-3} \text{ s}$ over five cardiac cycles of period $T = 0.95 \text{ s}$, by which point the flow rate and pressure fields have reached a periodic state. All results reported in this study are from the final cycle. 

\subsubsection{Point cloud and background mesh}
The boundary-fitted reference solution is obtained using the VMR model with the computational setup of Section~\ref{sec:aorta-setup}, employing the VMS formulation with weakly enforced no-slip conditions on the aortic surface. For the proposed framework, the nodes of the VMR surface mesh are sampled as a point cloud consisting of $72{,}882$ points, which is then immersed in a background domain as shown in Figure~\ref{fig:aorta-im-setup}. Imposing the inflow and outflow conditions requires care in this setting, since the inlet and outlet planes do not align with the background domain. While these conditions could be enforced weakly through surface integrals over the corresponding portions of the immersed boundary, this is undesirable for the inlet, where the velocity profile is typically prescribed strongly. Under weak enforcement, the inlet velocity remains part of the solution and satisfies the prescribed profile only approximately. This subjects the inlet to the same convective energy influx responsible for backflow divergence at open boundaries~\cite{Esmaily11compa, Bertoglio18Bench}. Moreover, the mass flux entering the domain would differ from that of the boundary-fitted reference, affecting the comparison. We therefore clip the background domain so that its boundary coincides with the desired inlet and outlet planes. For each opening, the corresponding points of the cloud are fitted to a plane, and a circle enclosing them is defined on that plane. Each circle is extruded outward along its normal, and the resulting cylinders are subtracted from the background domain using CAD boolean operations in Gmsh~\cite{Geuzaine09Gmsh}. The inlet and outlet planes are thereby exposed as fitted boundaries of the background mesh, on which the inflow and outflow conditions are imposed directly, while the aortic wall remains immersed and is treated with the proposed weak BC. 

The domain is then meshed with anisotropic near-wall refinement along the aortic surface, generating $4{,}360{,}005$ tetrahedral elements, of which only $1{,}851{,}268$ are active (intersected by or contained within the physical domain). The active element count is therefore lower than that of the boundary-fitted mesh, which consisted of $2.9$ million elements. The difference is notable because the boundary-fitted mesh was adaptively refined, which produced local refinement near the aortic arch and within the aorta, whereas the background mesh uses a uniform element size near the wall and in the interior. The element size and mesh statistics for the background mesh are reported in Table~\ref{tab:meshstatsaorta}.

\begin{table}[!t]\centering\small
    \centering
    \captionsetup{justification=centering}
    \newcommand{\tabincell}[2]{\begin{tabular}{@{}#1@{}}#2\end{tabular}}
    \caption{Background mesh statistics for patient-specific aorta. The element sizes are normalized by the largest diameter of the inlet cross-section ($D = 1.86$ cm).}
    \begin{tabular}{ccccr}
        \toprule
        \multicolumn{3}{c}{Near wall ($\times \, D$)}
        & \multirow{2}{*}{Domain ($\times \, D$)} 
        & \multirow{2}{*}{\tabincell{c}{Number of \\ active elements}} \\
        \cmidrule(lr){1-3}
        $h_\mathrm{NOR}$
        & $h_\mathrm{TAN}$
        & \tabincell{c}{BL height} 
        & & \\
        \midrule
        $0.01$  & $0.05$ & $0.05$  & $0.2$ & $1{,}851{,}268$  \\
        \bottomrule
    \end{tabular}
    \label{tab:meshstatsaorta}
\end{table}

\subsubsection{Results}
\begin{figure}[!t]
    \centering
    \begin{subfigure}{0.485\textwidth}\centering
        \includegraphics[width=\textwidth]{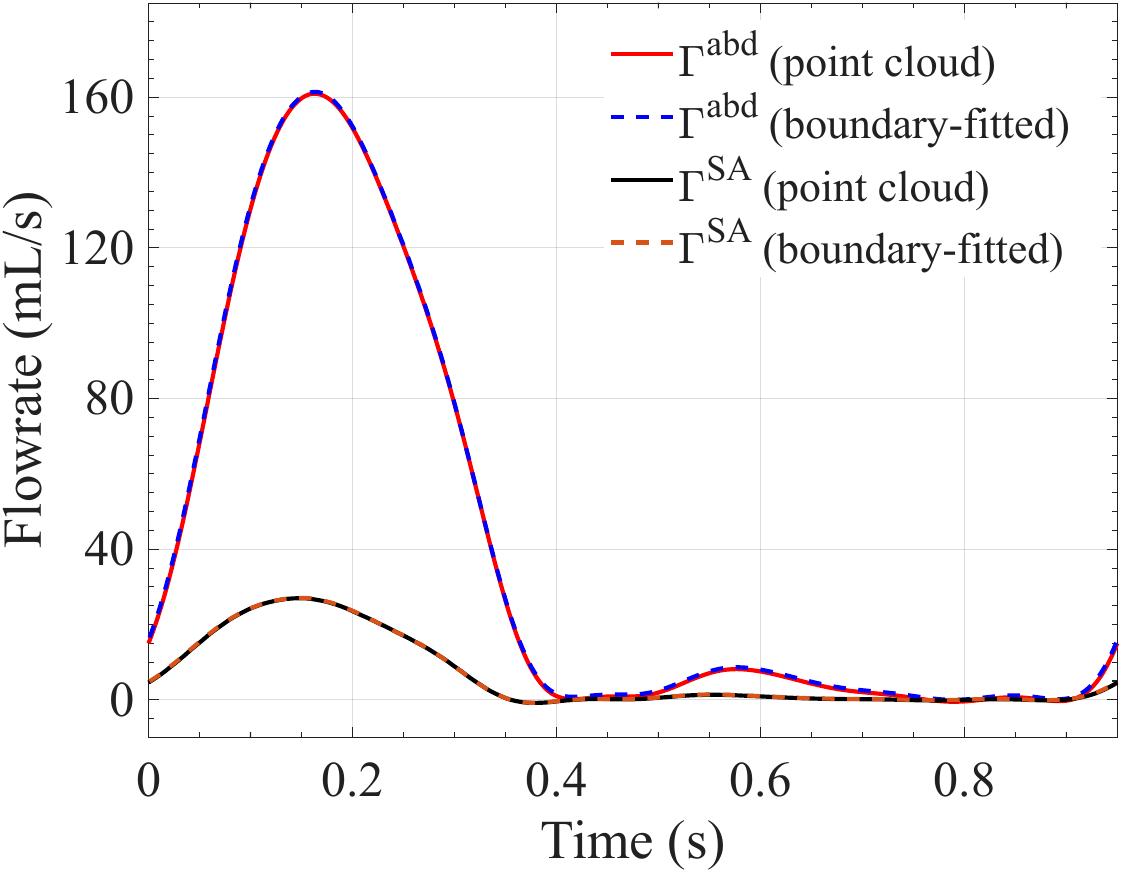}
        \caption{}
    \end{subfigure}
    \hspace{0.01\textwidth}
    \begin{subfigure}{0.485\textwidth}\centering
        \includegraphics[width=\textwidth]{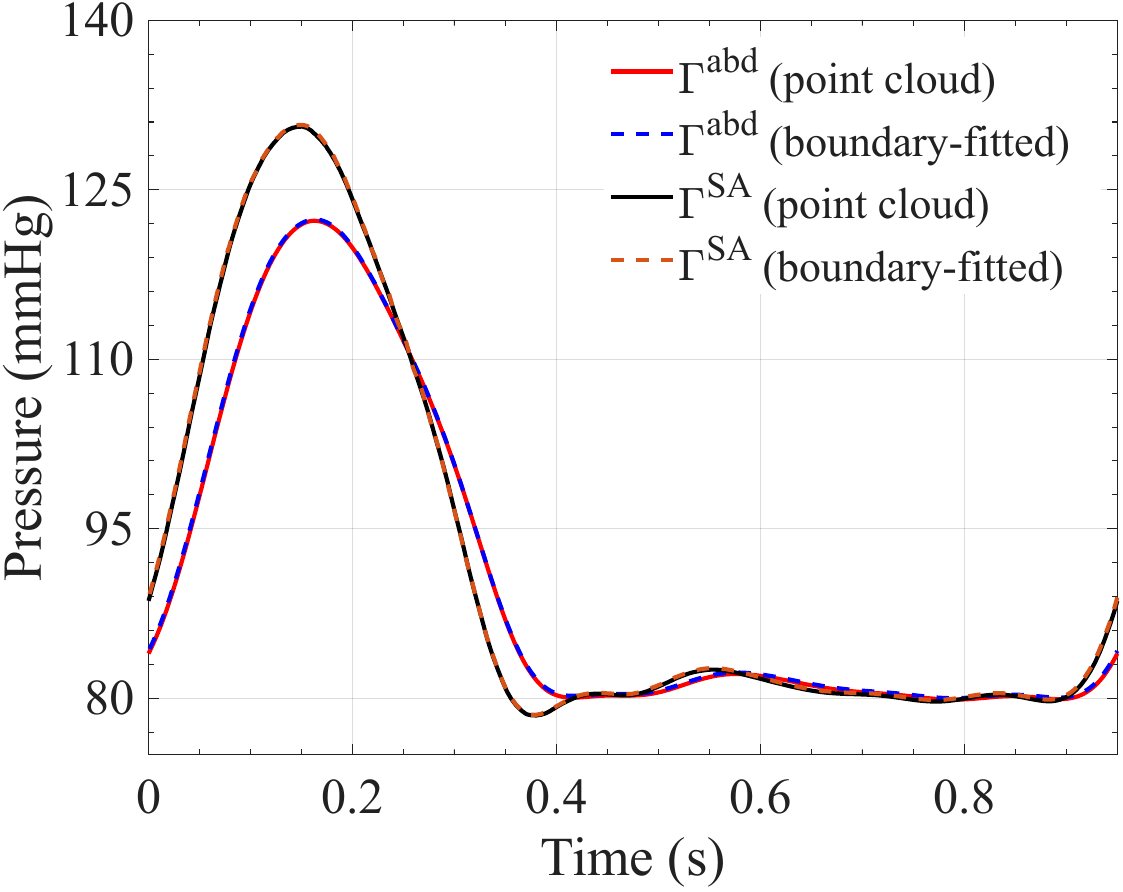}
        \caption{}
    \end{subfigure}
    \caption{(a) Flow rate and (b) area-averaged pressure on the supra-aortic and abdominal outlets over the final cardiac cycle compared with the boundary-fitted reference on patient-specific aorta.}
    \label{fig:aorta_outlets}
\end{figure}

\begin{figure}[!t]
    \centering
    \includegraphics[width=\textwidth]{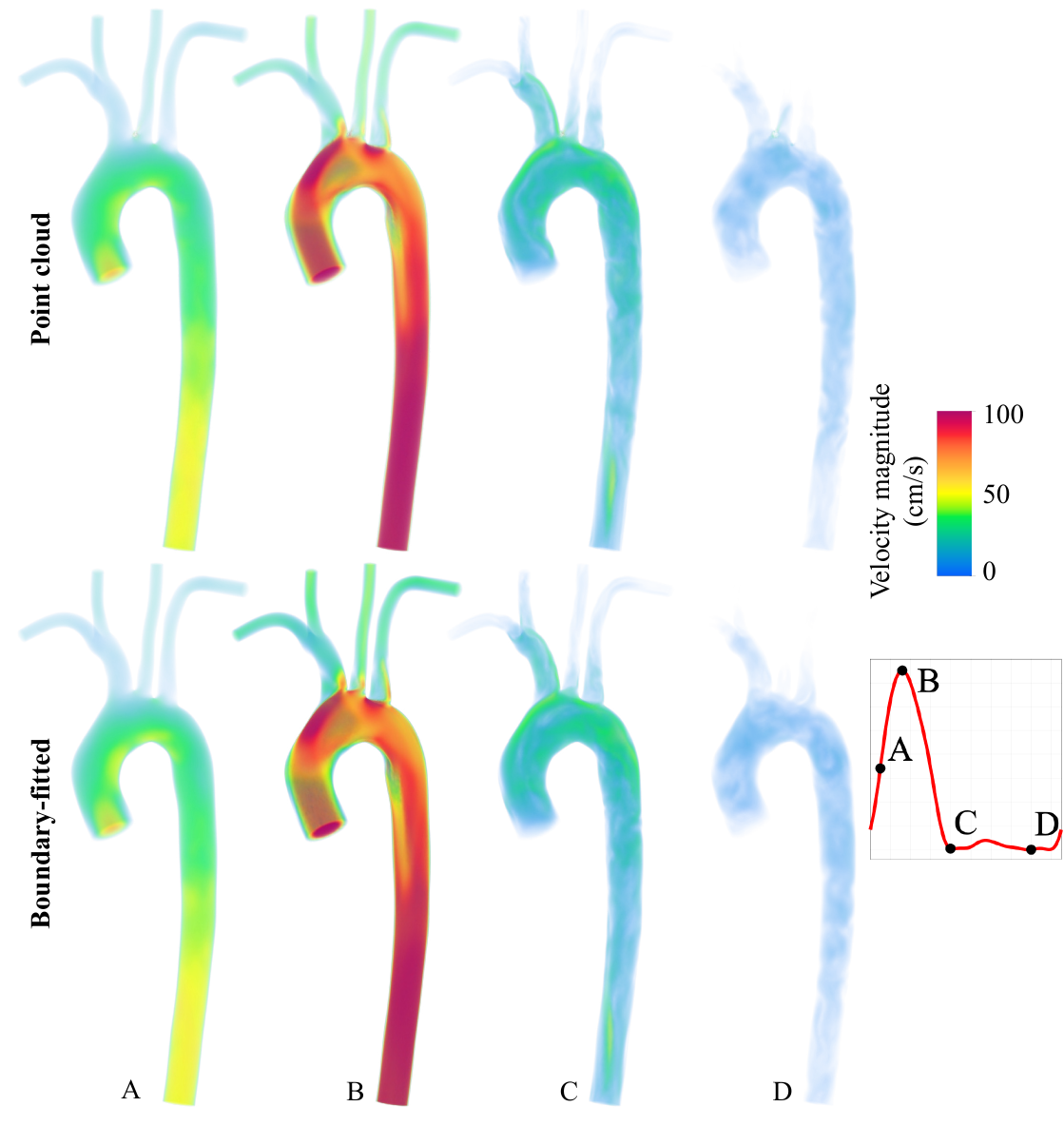}
    \caption{Volume contour of velocity magnitude at peak acceleration (A), peak systole (B), end of systole (C), and diastole (D) compared with boundary-fitted reference on patient-specific aorta. Top row: point cloud framework; bottom row: boundary-fitted reference.}
    \label{fig:aorta_vel}
\end{figure}

\begin{figure}[!t]
    \centering
    \includegraphics[width=\textwidth]{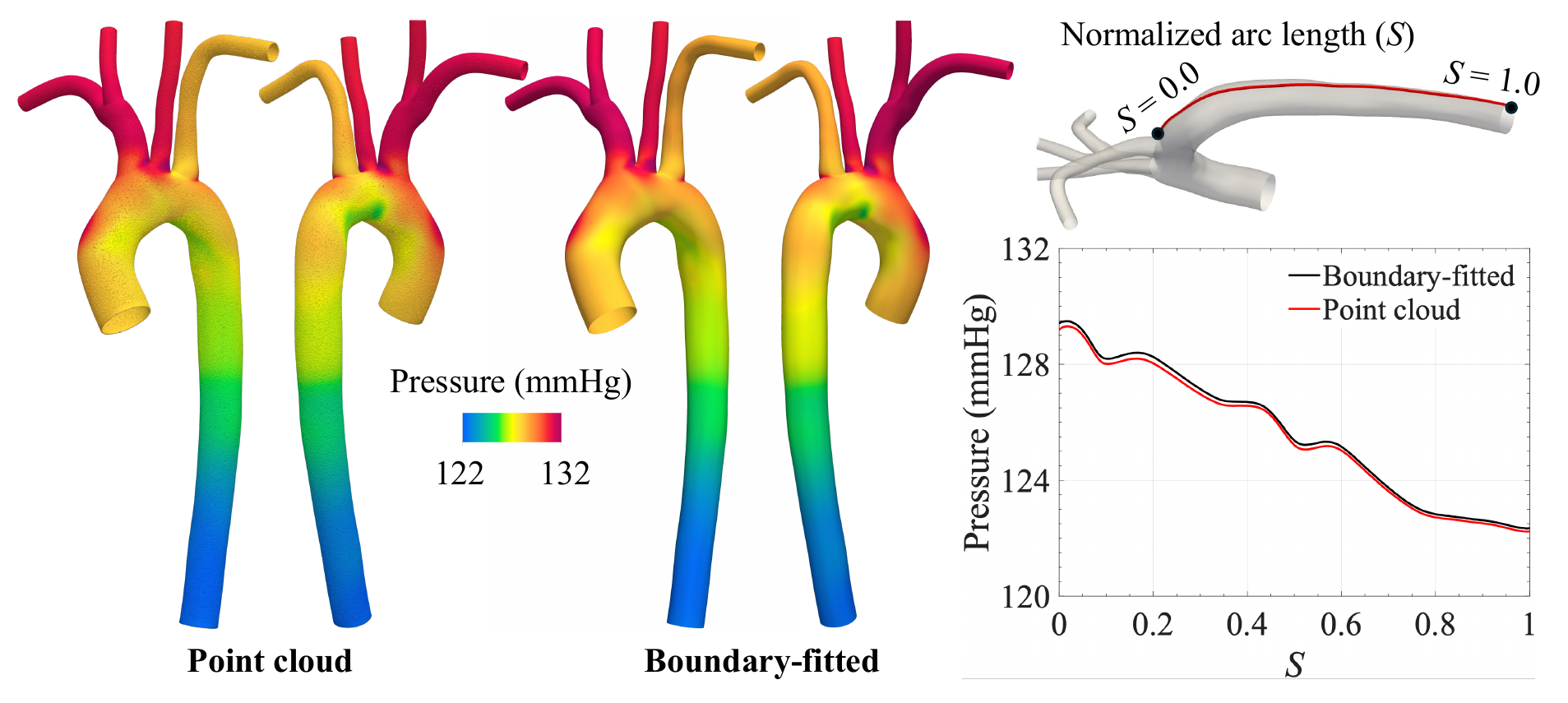}
    \caption{Instantaneous pressure distribution at peak systole evaluated on the point cloud compared with the boundary-fitted reference for patient-specific aorta. The plot shows the comparison of instantaneous pressure evaluated on the segment of outer curvature of the aorta, with the sampling curve indicated above.}
    \label{fig:aorta_pres}
\end{figure}

Figure~\ref{fig:aorta_outlets} compares the flow rate and the area-averaged pressure at the outlets over the final cardiac cycle against the boundary-fitted reference. Both quantities depict excellent agreement throughout the cycle, capturing the peak systolic ejection at $t\approx0.15$ s, the deceleration phase, and the diastolic period. The agreement in outlet flow rate confirms that the flow division between the supra-aortic arteries and abdominal aorta is reproduced by the point cloud framework. Moreover, the pressure agreement indicates that the resistance boundary conditions are enforced consistently in both settings. Figure~\ref{fig:aorta_vel} presents the volume contours of the velocity magnitude at four instants of the cardiac cycle, indicated on the inlet waveform: peak acceleration (A), peak systole (B), end of systole (C), and diastole (D). The point cloud and boundary-fitted simulations agree closely at all instants. At peak systole, both reproduce the high-velocity flow through the ascending aorta and the flow separation near the inner curvature. The pressure distribution on the aortic wall at peak systole is shown in Figure~\ref{fig:aorta_pres}. The point cloud solution reproduces the reference distribution over the entire surface, with elevated pressure in the ascending aorta and the supra-aortic branches and a monotonic decrease along the descending aorta. Figure~\ref{fig:aorta_pres} also shows a quantitative comparison of the pressure along the outer curvature of the aorta, parameterized by normalized arc length $S$. The two solutions agree closely along the entire curve, with mean and maximum differences of 0.15 and 0.23 mmHg, respectively. 

\begin{figure}[!t]
    \centering
    \includegraphics[width=\textwidth]{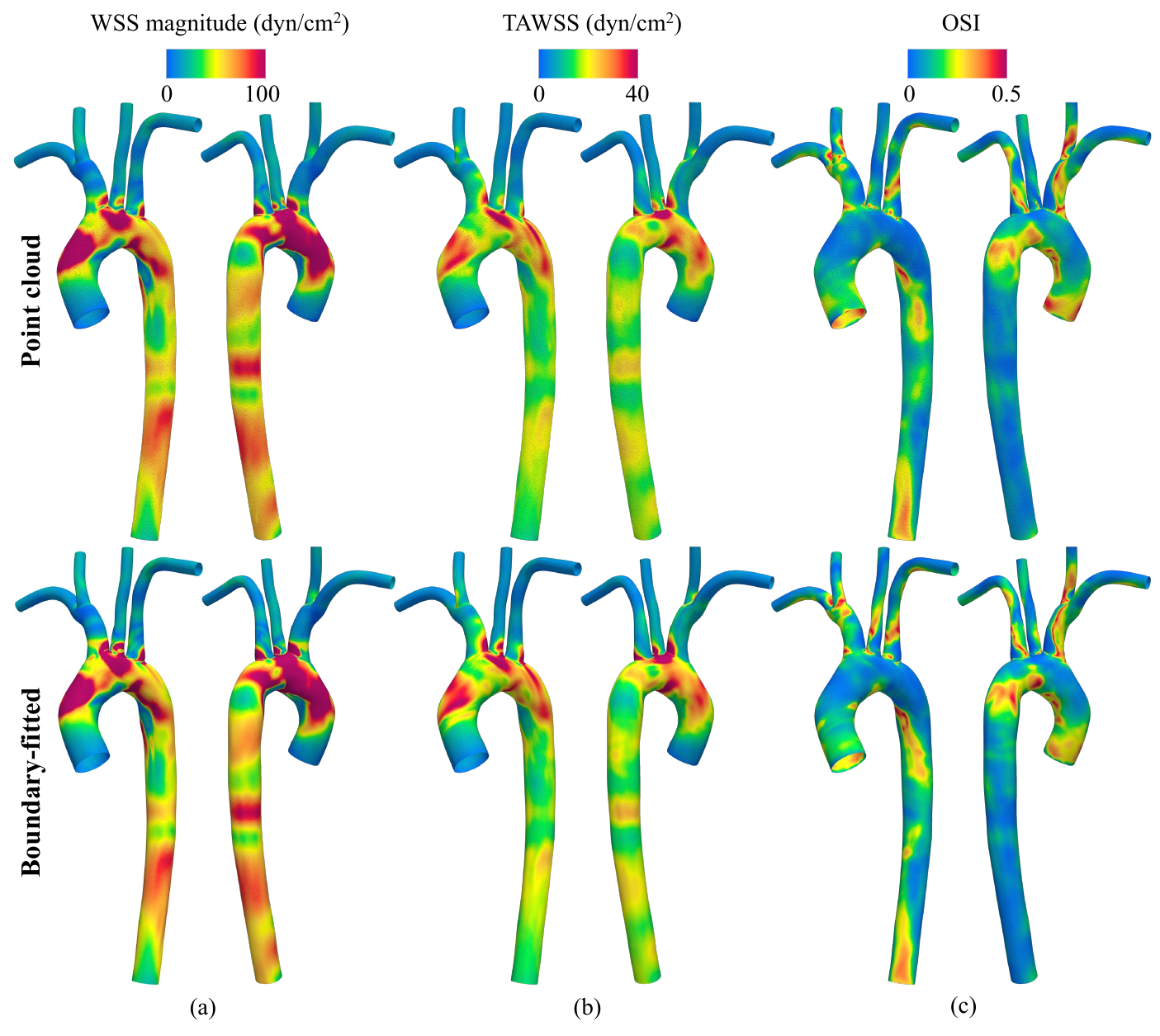}
    \caption{(a) Instantaneous WSS magnitude at peak systole, (b) TAWSS, and (c) OSI evaluated on the point cloud of the aortic wall (top row) and compared with the boundary-fitted reference (bottom row) for patient-specific aorta.}
    \label{fig:aorta_wss}
\end{figure}

The wall shear stress results are presented in Figure~\ref{fig:aorta_wss}, which compares the instantaneous WSS at peak systole, the time-averaged WSS (TAWSS), and the oscillatory shear index (OSI). The recovered WSS (Figure~\ref{fig:aorta_wss}a) is close to the reference field in both magnitude and spatial distribution, with elevated values near the arch and the abdominal aorta, and lower values near the supra-aortic arteries. The same agreement is observed for TAWSS (Figure~\ref{fig:aorta_wss}b), including the location and extent of the low-TAWSS regions that are of primary clinical interest. The OSI distribution (Figure~\ref{fig:aorta_wss}c) shows overall good agreement as well, with slight local variations near the inner curvature. This is the most demanding of the three comparisons, as OSI depends on the directional history of the WSS over the entire cycle and is therefore sensitive to any spatial or temporal variation in the underlying traction. The minor differences in the results are most likely attributable to the differing discretizations. The boundary-fitted mesh was adaptively refined, whereas the background mesh uses a uniform element size at the wall and in the interior, without local adaptation. Even with this simpler discretization and with fewer degrees of freedom, the point cloud framework closely matches the reference WSS, TAWSS, and OSI distributions, including low-TAWSS and elevated-OSI regions of primary clinical interest. 

\section{Conclusions}
\label{sec:conclusion}
In this paper, we presented a point cloud-based CFD framework for performing flow analysis directly on geometries represented by point clouds, with emphasis on the accurate recovery of wall shear stress. The framework is based on immersogeometric analysis using the VMS formulation for the incompressible Navier--Stokes equations, ghost penalty stabilization for small-cut elements, non-symmetric Nitsche-based weak BC, and a stress recovery procedure for the WSS. The non-symmetric weak BC, with near-wall modeling for parameter estimation, was presented in detail and compared against its symmetric variant in the immersed setting. The WSS recovery method was developed in full, comprising a geometry-following element-based patch construction and a least-squares formulation that constrains the recovered stress field with variationally consistent traction supplied by the weak BC. The framework was validated on canonical benchmark cases of Hagen--Poiseuille pipe flow and laminar flow around a sphere, showcasing excellent mesh convergence and agreement with analytical and boundary-fitted results. It was then applied to turbulent flow past a sphere at $Re = 3700$, where the recovered skin friction agrees with published reference data. Finally, it was used to simulate a patient-specific aorta, where the velocity, pressure, WSS, TAWSS, and OSI distributions agree closely with boundary-fitted computation on the same geometry. 

The proposed method now enables point cloud-based CFD not only to perform direct flow analysis on scanned real-world geometries, without surface reconstruction, geometry cleanup, or boundary-fitted mesh generation, but also to compute accurate near-wall quantities, particularly wall shear stress. This opens the way to future extensions of the framework, including multiphysics applications such as fluid--structure interaction, and to the development of digital twins built on real-world scans rather than on idealized models.

\section*{Declaration of Competing Interest}
The authors declare that they have no known competing financial interests or personal relationships that could have appeared to influence the work reported in this paper.

\section*{Acknowledgments}
This work was supported in part by the National Heart, Lung, and Blood Institute of the National Institutes of Health under award number R01HL184128 and by the National Science Foundation under award number DMS-2436623. This support is gratefully acknowledged. We also thank the Texas Advanced Computing Center (TACC) at The University of Texas at Austin for providing computational resources that have contributed to the research results reported within this paper.

\appendix
\section{Effect of quadrature on internal flows}\label{App:AppendixA}
\setcounter{figure}{0}
\setcounter{table}{0}

The cut quadrature of Section~\ref{sec:cut-quad} was motivated by efficiency. This appendix shows that it also removes an accuracy limitation of adaptive subdivision quadrature that becomes significant in internal flows driven by pressure boundary conditions. 

In the variational form~\eqref{eq:weak_form}, the pressure enters the momentum residual of the VMS formulation through volume and boundary terms. Let $(\cdot \, , \cdot)_{\Omega_\text{phys}}$, $(\cdot \, , \cdot)_{\Gamma^b}$, and $(\cdot \, , \cdot)_{\Gamma^\text{N}}$  denote the exact $L^2$ inner products over the physical domain, the immersed boundary, and the Neumann boundary, respectively, and let $(\cdot \, , \cdot)_{\Omega_\text{phys}}^\mathcal{Q}$, $(\cdot \, , \cdot)_{\Gamma^b}^\mathcal{Q}$, and $(\cdot \, , \cdot)_{\Gamma^\text{N}}^\mathcal{Q}$ denote their evaluation by the integration scheme $\mathcal{Q}$, such that
\begin{align}
(f, g)^\mathcal{Q}_\mathcal{W} = \sum_k w_k^\mathcal{W} f(\mathbf{x}_k^\mathcal{W}) \, g(\mathbf{x}_k^\mathcal{W}) \text{ ,} \quad \mathcal{W} \in \{\Omega_\text{phys}, \Gamma^b, \Gamma^\text{N} \} \text{ ,}
\end{align}
where $\mathbf{x}_k^\mathcal{W} \in \mathcal{W} $ and $w_k^\mathcal{W}$ are quadrature points and weights generated by $\mathcal{Q}$ on $\mathcal{W}$. The pressure contribution to the discrete momentum residual is given as
\begin{align}
r_{p}^\mathcal{Q} (\mathbf{w}^h; p^h, p_0) =  
(p^h , \nabla \cdot \mathbf{w}^h )_{\Omega_\text{phys}}^\mathcal{Q}
- (p^h  , \mathbf{w}^h \cdot \mathbf{n})_{\Gamma^b}^\mathcal{Q}
- (p_0  , \mathbf{w}^h \cdot \mathbf{n})_{\Gamma^\text{N}}^\mathcal{Q} \text{ ,}
\end{align}
where $p_0$ is the prescribed pressure at the Neumann boundary $\Gamma^\text{N}$. Consider a discrete solution $(\mathbf{u}^h, p^h)$ and shift the pressure by a constant $\bar{p}$, both in the Neumann condition and in the solution. Since a uniform pressure has zero gradient, the physical problem is unchanged with $(\mathbf{u}^h, p^h + \bar{p})$ as a solution. However, by linearity, the difference in residual is 
\begin{align}
   \label{eq:geom-error}
   r_{p}^\mathcal{Q} (\mathbf{w}^h; p^h + \bar{p}, p_0 + \bar{p})  - r_{p}^\mathcal{Q} (\mathbf{w}^h; p^h, p_0) = \bar{p} \, \epsilon_\mathcal{Q}(\mathbf{w}^h) \, 
\end{align}
where
\begin{align}
    \epsilon_\mathcal{Q}(\mathbf{w}^h) = (1 , \nabla \cdot \mathbf{w}^h )_{\Omega_\text{phys}}^\mathcal{Q}
- (1  , \mathbf{w}^h \cdot \mathbf{n})_{\Gamma^b}^\mathcal{Q}
- (1  , \mathbf{w}^h \cdot \mathbf{n})_{\Gamma^\text{N}}^\mathcal{Q}\text{ }
\end{align}
defines a geometric error, which only depends on the geometry and the quadrature rule considered. With exact integration, $\epsilon_\mathcal{Q}(\mathbf{w}^h)  = 0$ by the divergence theorem. However, when the volume quadrature (adaptive subdivision of the cut element) and the surface quadrature (the point cloud with its areas or from a surface mesh) are constructed independently, the two quadrature errors do not cancel, and $\epsilon_\mathcal{Q}$ is non-zero. The shifted pressure therefore leaves a spurious momentum residual, concentrated on the cut elements and growing linearly with $\bar{p}$, which propagates as a nonphysical mass flux through the immersed boundary. In external flows, the relevant pressure level is of the same order as the dynamic pressure and the error is small. However, in internal flows with prescribed pressure conditions, $\bar{p}$ can exceed the dynamic pressure by orders of magnitude, amplifying the geometric error. The cut quadrature avoids the problem by construction. Since the surface and volume rules derive from the same sub-cell decomposition, the surface rule is exactly the boundary of the volume rule and the divergence theorem holds to the accuracy of the local geometric approximation. 

\begin{figure}[!t]
    \centering
    \begin{subfigure}{0.475\textwidth}\centering
        \includegraphics[width=\textwidth]{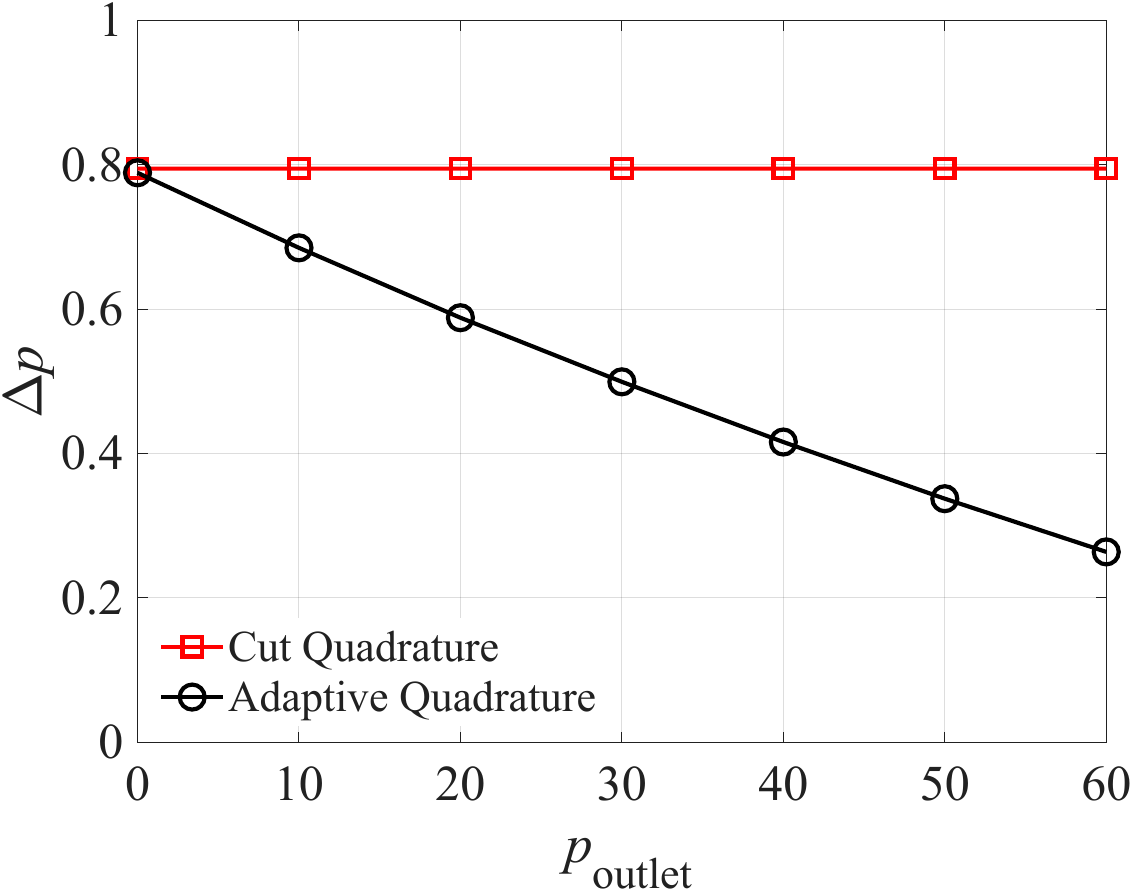}
        \caption{}
    \end{subfigure}
    \hspace{0.01\textwidth}
    \begin{subfigure}{0.485\textwidth}\centering
        \includegraphics[width=\textwidth]{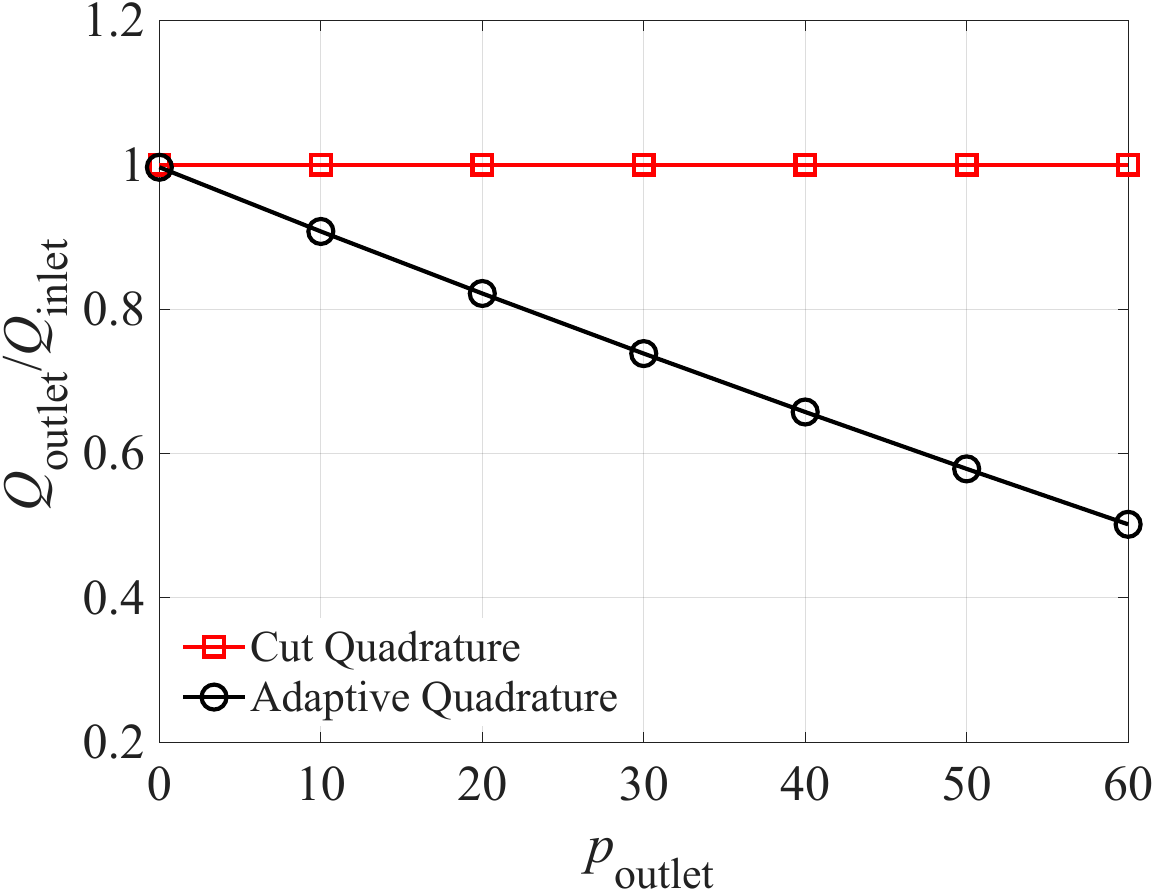}
        \caption{}
    \end{subfigure}

    \caption{(a) Pressure drop $\Delta p$ and (b) outlet-to-inlet flow rate ratio as functions of the prescribed outlet pressure $p_\text{outlet}$, for cut quadrature and adaptive subdivision quadrature. The analytical values, $\Delta p = 0.8$ and $Q_\text{outlet}/Q_\text{inlet} = 1$, are independent of $p_\text{outlet}$.}
    \label{fig:pipe-CQvsAQ-plot}
\end{figure}

\begin{figure}[!t]
    \centering
    \includegraphics[width=\textwidth]{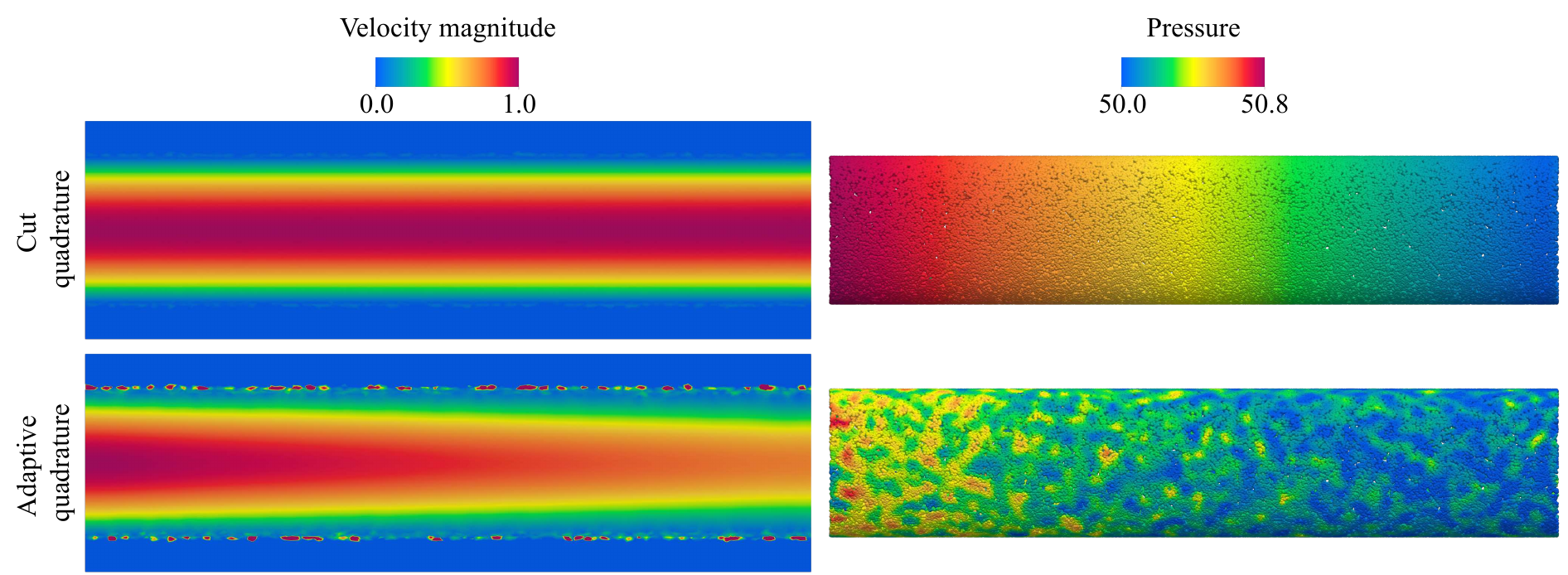}
    \caption{Velocity magnitude on the central cross-section (left) and pressure evaluated on the point cloud (right) at $p_\text{outlet} = 50$ on mesh IM1 for cut quadrature (top) and adaptive quadrature (bottom) with subdivision level of 2. }
    \label{fig:pipe-CQvsAQ}
\end{figure}

We demonstrate this with a patch test on the Hagen--Poiseuille flow of Section~\ref{sec:hagen}, using mesh IM1. All parameters are held fixed while the traction-free BC is replaced with an outlet pressure increasing from 0 to 60, which shifts the pressure without changing the physical problem. The analytical pressure drop remains $\Delta p = 0.8$ and mass conservation requires $Q_\text{outlet}/Q_\text{inlet} = 1$ for every $p_\text{outlet}$. The cut quadrature rule is compared with the adaptive quadrature with 2 subdivision levels. Figure~\ref{fig:pipe-CQvsAQ-plot} shows $\Delta p$ and $Q_\text{outlet}/Q_\text{inlet}$ for different $p_\text{outlet}$. With cut quadrature, both quantities remain at their analytical values over the entire range, while, with adaptive quadrature, they deviate significantly as $p_\text{outlet}$ increases. Moreover, a linear relationship between the error for both quantities and $p_\text{outlet}$ is observed, consistent with Eq.~\eqref{eq:geom-error}. Figure~\ref{fig:pipe-CQvsAQ} demonstrates the velocity and pressure field obtained with both quadrature rules at $p_\text{outlet}=50$. The cut quadrature solution recovers the parabolic profile and a smooth linear pressure drop along the wall. In contrast, the adaptive quadrature solution exhibits spurious velocity at the immersed boundary, velocity decay along the pipe as mass leaks through the wall, and wall pressure dominated by noise. Within the adaptive framework, this error can only be reduced, at increasing cost, by simultaneously refining both volume and surface quadrature. The cut quadrature eliminates it by construction. This is of direct relevance to the aorta case of Section~\ref{sec:aorta}, where the physiological pressure level is two orders of magnitude larger than the pressure variation along the vessel.

\small
\bibliographystyle{vancouver-mch}
\bibliography{refs-mch}

\end{document}